\documentclass[11pt]{article}

\usepackage[authoryear]{natbib}

\usepackage{times}
\usepackage{graphicx} % Required for inserting images
\usepackage{mathtools}
\usepackage{amsmath}
\usepackage{amssymb}
\usepackage{float}
\usepackage{longtable}

\usepackage{tocbasic}
\usepackage{caption}

\usepackage{subcaption}
\usepackage{bm}
\usepackage[margin=1.2in]{geometry}
\usepackage{hyperref}
\usepackage{lscape}
\usepackage{multirow}
\usepackage{titlesec}
\usepackage{siunitx}
\usepackage[USenglish]{babel}
\usepackage{algorithm} 
\usepackage{algpseudocode} 
\usepackage[intoc]{nomencl}
\usepackage{hyphenat}
\usepackage{authblk}

\newcommand{\bx}{\bm{x}}

\newcommand{\rmd}{\mathrm{d}}

\newcommand{\balpha}{\bm{\alpha}}
\newcommand{\bvarphi}{\bm{\varphi}}

\newcommand{\R}{\mathbb{R}}
\newcommand{\E}{\mathbb{E}}
\newcommand{\calX}{\mathcal{X}}
\newcommand{\calY}{\mathcal{Y}}
\newcommand{\D}{\mathcal{D}}
\renewcommand{\O}{\mathcal{O}}
\newcommand{\GP}{\mathcal{GP}}
\newcommand{\calD}{\mathcal{D}}
\newcommand{\Cov}{\mathbb{C}\text{ov}}

\newcommand{\iain}{\textcolor{black}}
\newcommand{\sylvain}{\textcolor{black}}
\newcommand{\denis}{\textcolor{black}}

\begin{document}

\date{}
\title{PFEM-GP-dPHS : a finite element framework for combining Gaussian processes and infinite-dimensional port-Hamiltonian systems}

% Authors with different addresses:
\author[a]{Florian Courteville} 
\author[b]{Iain Henderson}
\author[b]{Denis Matignon}
\author[c]{Sylvain Dubreuil}

\affil[a]{KTH royal institute of technology, Stockholm, Sweden \\ Institut Laue-Langevin, 38000, Grenoble, France \texttt{fbco@kth.se}}

\affil[b]{Fédération ENAC ISAE-SUPAERO ONERA, Université de Toulouse, 31000, Toulouse, France \texttt{\{iain.henderson, denis.matignon\}@isae-supaero.fr}}

\affil[c]{DTIS, ONERA, Université de Toulouse, 31000, Toulouse, France \\ Fédération ENAC ISAE-SUPAERO ONERA, Université de Toulouse, 31000, Toulouse, France \texttt{sylvain.dubreuil@onera.fr}}

\maketitle

\begin{abstract}%
In order to learn distributed port-Hamiltonian systems (dpHs) using Gaussian processes (GPs), the partitioned finite element method (PFEM) is combined with the Gp-dpHs method.
By following a late lumping approach, the discretization of the functional hyperparameters of the GP prior over the Hamiltonian functional is chosen independently from the discretization of the dpHs, thus reducing the numerical complexity of our method. We next model the mean of the GP prior of the Hamiltonian as a quadratic form, enabling the GP kernel to focus on the nonlinear part of a given dpHs. We illustrate our method on a nonlinear one dimensional wave equation with unknown physical parameters (tension and linear mass). % and show that its physical spatially dependent parameters (tension and linear mass) can be estimated using our method.
\end{abstract}
\textbf{Code} : \href{https://github.com/olimarlouis/PFEM-GP-dpHs.git}{https://github.com/olimarlouis/PFEM-GP-dpHs.git}  \\
\textbf{Keywords :}  Physics-informed Machine Learning, Distributed Parameter Port-Hamiltonian System, Partitioned Finite element Method, System Identification%

\section{Introduction \iain{and contribution}}

%Introduire rapidmeent GP-PHS, puis vouloir l'appliquer a dim infinie discrétiser -> comment faire? PFEM puis hamilton puis GP-PHS? TROp de paramètre, donc GP sur la fonctionnel = paramètre def continue = plein de base

%Pour le context, St Venant, ... a trouver 

%\sylvain{On utilise beaucoup "we will ...", on m'a toujours dit de ne pas écrire au futur dans une publi, un avis sur la question ? (ça ne me dérange pas de laisser comme ça)}

The modelling and identification of dynamical systems is a crucial task in a broad range of domains such as physics, engineering and applied mathematics. \iain{The mathematical framework of Port-Hamiltonian Systems (PHS) is especially well suited to describe interacting dynamical systems. Gaussian processes are functional Bayesian models that can be used for both model approximation and identification while also providing uncertainty quantification \citep{RW}.}
In this context, T. Beckers \iain{and co-authors} introduced Gaussian process port-Hamiltonian systems (GP-PHS) for ordinary differential equations (ODEs) \citep{Beckers}. This method couples Gaussian processes (GPs) with the PHS formalism to obtain physics-based GP priors, so that the result obtained from the data-driven modelling respects the \iain{PHS structure} of the unknown system, which has been proved useful in \citep{Mahboubi_2025}. Attempts to extend this approach to \iain{partial differential equations (PDEs), known as }distributed port-Hamiltonian system (dPHS) have been proposed in \citep{tan24,NAPI-MPC,Pnp-PIMP}, leading to the Gaussian process distributed port-Hamiltonian system (GP-dPHS) method. %\iain{faire une coupure ici pour une section contributions?} 
\sylvain{These two contributions apply the GP-PHS approach to the discretised system obtained using an ad hoc (finite difference) discretization method leading to an approximate finite-dimensional system (an ODE system). 
}

\iain{In the present article, we make use of PFEM \citep{PFEM}, a finite element method  designed to respect the structure of a given dPHS at the discrete level. However, applying the GP-PHS method directly to the resulting ODE system would result in a prohibitively large number of hyperparameters. To solve this issue, we adapt the physics-based priors from the GP-dPHS method within the discretization method, PFEM, to obtain a prior whose hyperparameters are defined continuously, and can then be discretized with their own adapted method.}   % to the discretization method employed is chosen to ensure that the finite-dimensional system respects some properties of the original dPHS. 

%Second, applying the GP-PHS method to the finite-dimensional system may lead to a number of parameters to identify that is much higher that what would be necessary to describe the regressed system. In the present contribution it is proposed, first, to rely on the partitioned finite element method (PFEM) \citep{PFEM} that preserves important properties such as the power balance of the finite dimensional approximate system. And second to adapt the physics-based priors from the GP-dPHS method within the discretization method, PFEM, to obtain a prior whose parameters are defined continuously, and can then be discretized with their own adapted method. Hence this contribution explores the "model by GP, then discretize" paradigm rather than the "discretize, then model by GP" paradigm. 
%} However when dealing with PDE and dPHS, discretization methods need to be used to obtain a approximate finite-dimensional system (an ODE system) that can by solved numerically. One such efficient method for dPHS is the partitioned finite element method (PFEM) \citep{PFEM} that preserves PHS \sylvain{properties} \iain{on the finite-dimensional approximate system}, such as the power balance. 
%

\section{Background}
%\iain{Direct introduire des PHS de dim infinie uniquement.}
%\iain{Courte section GP : indexé par des fonctions, puis noyau à val opérateurs donner ref.}
We first introduce GPs and PHS, and then present the GP-dPHS method.
%Introduire les dPHS, puis les GP au noyau mat/operateur. Ensuite on va introduire le GP-dPHS avec la fonctionnel sous GP, et loperateur JR replace par approximation. Efin lors de PFEM, la fonctionnel se faire discretiser avec les consitutive, et JR se fait discretiser de l'autre coté aussii, puis on assemble.

% \subsection{Port-Hamiltonian systems}

%  The dynamics of a PHS are described by:

% \begin{equation}
%     \begin{cases}
%         \dbx(t) = [J(\bx(t))-R(\bx(t))]\nabla_xH(\bx(t)) +G(\bx(t))\bu(t) \\ 
%         \by(t)=G^{\top}(\bx(t))\nabla_xH(\bx(t))
%     \end{cases}
% \end{equation}

% Above, $\bx\in\R^n $ is the state (or energy variable) of the system and $\dbx := \frac{d \bx}{dt}$. $\nabla_xH(\bx)\in\R^n$ is called the co-energy variable. $\bu\in\R^m$ and $\by\in\R^m$ are the input and output of the system. $J : \R^n\rightarrow \R^{n\times n}$ is a skew-symmetric matrix which describes the exchange of energy with the storage components. $R : \R^n\rightarrow \R^{n\times n}$ is a symmetric non-negative matrix which describes the exchange with the dissipative components. $G : \R^n\rightarrow \R^{n\times m}$ is a rectangular matrix which describes the exchange with the environment. $H : \R^n \rightarrow\R$ is a smooth Hamiltonian function in $\mathcal{C}^\infty$.

\subsection{Scalar, vector and function-valued Gaussian processes}\label{sec:gp}

In this section we introduce scalar, vector and function-valued Gaussian processes indexed by a Hilbert space, and the properties that we will use.
 \iain{Let $\mathcal{X}$ be a Hilbert space, and $\mathcal{Y}$ be either $\R^n, n\geq 1$ or a Hilbert space}. Let ($\Omega$, $\mathcal{F}$, $\mathbb{P}$) be a probability space and \iain{$f = (f(\balpha))_{\balpha \in \mathcal{X}}$ a family of $\calY$-valued random variables, $f(\balpha) : \Omega \rightarrow \calY$}. \iain{Below, $f $ will be a  Gaussian process; its distribution will be determined by its mean function $m : \calX\rightarrow \calY, \ \alpha \mapsto \mathbb{E}(f(\balpha))$ and covariance function (or kernel) $k$, per detailed below in the cases where $\calY = \R,\  \calY = \R^n$ and $\calY$ is a Hilbert space. We summarize this as $f\sim \GP (m,k)$. }

\paragraph{Scalar-valued GP \citep{RW}}
If $\calY = \R$ then $f \sim \GP (m,k) $ is a scalar-valued GP, meaning that for all $\balpha_1, \ldots, \balpha_N\in \calX$, $(f(\balpha_1), \ldots, f(\balpha_N))$ follows a multivariate Gaussian distribution. The kernel $k$ is given by $k: \calX\times\calX \rightarrow \R, (\balpha,\balpha') \mapsto \Cov(f(\balpha),f(\balpha'))$. For $k$ to be the covariance function of a GP, it is necessary and sufficient that $k$ be symmetric and positive semi-definite (p.s.d.), i.e. %for all $N\in \mathbb{N}$, $\balpha_1, \ldots,\balpha_N \in \mathcal{X}$ and $c_1,\ldots,c_N\in \R$, we have
\begin{equation}
    \forall N\in \mathbb{N}, \ \ \ \forall \balpha_1, \ldots,\balpha_N \in \mathcal{X},\ \ \ \forall \ c_1,\ldots,c_N\in \R, \ \ \ \ \ \sum_{i=1}^N\sum_{j=1}^Nc_i c_jk(\balpha_i,\balpha_j) \geq 0.
    \label{eq:PSD}
\end{equation} 
\iain{GPs can be conditioned on measurement data.} For this,
%\iain{The GP $f$ can be conditioned as follow \citep{Roustant}.} 
consider outputs $y_i$ \iain{at location $\balpha_i$}, affected by Gaussian noise $\eta_i \sim \mathcal{N}(0,\sigma^2)$, with $(n_i)$ i.i.d. : $y_i = f(\balpha_i) +\eta_i$. This yields a training set $\mathcal{D} = \{X,Y\}$ of input/output data, $X = [\balpha_1,\ldots \balpha_N]$ and $Y = [y_1,\ldots y_N]$. Then $f \sim \GP(m,k)$ can be conditioned on $\calD$, yielding a \textit{posterior} GP with explicit mean and covariance functions $m_\D(\alpha)=\E(f(\balpha)|\mathcal{D})$ and $k_\calD(\balpha,\balpha') = \Cov(f(\balpha),f(\balpha')|\calD)$ (see \citep{Roustant}, Chapter 4 for further details).  %, we obtain the so-called Kriging posterior
%\begin{equation}
%    \E(f(\balpha)|\mathcal{D}) = m(\balpha^*) + k(\balpha^*,X)^\top K^{-1}Y_0,
%\end{equation}

%\begin{equation}
%   \Cov(f|\balpha^*,{\balpha^*}',\mathcal{D}) = k(\balpha^*,{\balpha^*}') - k(\balpha^*,X)^\top K^{-1}k({\balpha^*}',X).
%\end{equation}

%Above $K$ is the Gram matrix whose elements are given by $K_{i,j} =k(X_{i,:},X_{j,:})+\delta(i,j)\sigma²$ where $\delta$ is the Kronecker delta. $Y_0$ is the mean-adjusted output data $Y_0 = [y_1-m(\balpha_1),\ldots y_N-m(\balpha_N)]$.

\paragraph{Vector-valued GP \citep{matrix_kernel}}

If $\calY = \R^n$ then $f \sim \GP (m,k) $ is a vector-valued GP with \textit{vector-valued} mean function $m : \mathcal{X} \rightarrow \R^n$ and symmetric \textit{matrix-valued} kernel $k : \calX \times \calX \rightarrow  \R^{n\times n}$; the entries of $k(\balpha,\balpha')$ describe the covariance between the coordinates of $f(\balpha)$ and $f(\balpha')$, through $k(\balpha,\balpha')_{ij} = \Cov(f(\balpha)_i,f(\balpha')_j)$. The p.s.d. property of $k$ reads as
\begin{equation}
    \forall N\in \mathbb{N}, \ \ \ \forall \balpha_1, \ldots,\balpha_N \in \mathcal{X},\ \ \ \forall \ c_1,\ldots,c_N\in \R^n, \ \ \ \ \ \sum_{i=1}^N\sum_{j=1}^Nc_i^\top k(\balpha_i,\balpha_j)c_j \geq 0.
    \label{eq:PSD-matrix}
\end{equation} 

%Then with measures affected by Gaussian noise $\eta_i \sim \mathcal{N}(0,\sigma^2 I_n)$, with $(n_i)$ i.i.d., the Gram matrix for conditioning became $K\in\R^{nN\times nN} $  whose sub-matrices are given by $K_{in:(i+1)n,jn:(j+1)n} =k(X_{i,:},X_{j,:})+\delta(i,j)\sigma^2 I_n$, and $Y_0\in\R^{nN}$ sub-vector are given by $(Y_0)_{in:(i+1)n} = y_i-m(\balpha_i)$ .

\paragraph{Function-valued GP \citep[Appendix A]{BATLLE2024}}

Finally, if $\calY$ is a Hilbert space then $f \sim \GP (m,k) $ is a $\calY$-valued GP with function mean function $m : \mathcal{X} \rightarrow \calY$ and symmetric kernel $k : \mathcal{X} \times \calX \rightarrow  \mathcal{L}(\calY)$. The p.s.d property of $k$ reads as %is that for all $N\in \mathbb{N}$, $\balpha_1, \ldots,\balpha_N \in \mathcal{X}$ and $c_1,\ldots,c_N\in H$, we have using the scalar product of $H$
\begin{equation}
 \forall N\in \mathbb{N}, \ \ \ \forall \balpha_1, \ldots,\balpha_N \in \mathcal{X},\ \ \ \forall \ c_1,\ldots,c_N\in \calY, \ \ \ \ \ \sum_{i=1}^N\sum_{j=1}^N\langle c_i, k(\balpha_i,\balpha_j)c_j\rangle_\calY \geq 0.
    \label{eq:PSD_op}
\end{equation}

\subsection{Distributed Port-Hamiltonian Systems (dPHS)}

Adding input/output ports to a Hamiltonian system leads to a port-Hamiltonian system. The port describes the exchange of energy between the energy storage components (spring, kinetic energy...), the dissipative components (friction, resistor...) and the environment. If the energy variables $\balpha \in \mathcal{X}$ are distributed within a continuous space domain $\mathcal{O}$, so that $\balpha :x \in \mathcal{O} \mapsto\balpha(x)$ are function of space, then we talk about distributed port-Hamiltonian system \citep{Vdshaft_infinite}. To write down the dPHS equation of a system, we first need to define a Hamiltonian functional $ \mathcal{H}$ describing the total energy of the system, with

\begin{equation}
    \mathcal{H}(\balpha) = \int_{\O} H(\bx,\balpha) d\bx
\end{equation}

where $H : \O \times \mathcal{X} \rightarrow \R$ is the energy density. Then, using the variational derivative of $\mathcal{H}$, we define the co-energy variables $\bm{e} := \delta_{\balpha}\mathcal{H}$, see \citep{Vdshaft_infinite}. We call constitutive relation the relation between $\bm{e}$ and $\balpha$. Finally, the dPHS form is obtained by rewriting the PDE of the system using the energy and co-energy variable such that

\begin{equation}\label{eq:dphs-balance}
    \partial_t\balpha := \dot{\balpha} = (\mathcal{J}-\mathcal{R})\bm{e}(\balpha),   \qquad
    \frac{d}{dt}\mathcal{H}(\balpha) = - \int_{\cal O} \bm{e}(\balpha).\mathcal{R}\bm{e}(\balpha) + \langle u_\partial, y_\partial \rangle_{\partial \cal O}.
\end{equation}

Above $\mathcal{J}$ is a \denis{constant skew-symmetric differential operator} which describes the exchange of energy with the storage components and $\mathcal{R}$ is a \denis{constant nonnegative symmetric operator} which describes the exchange with the dissipative components. \denis{Collocated input-output ($u_\partial$, $y_\partial$) behaviour} describing the interaction with the environment are defined by the boundary conditions of the problem (such as Dirichlet control, Neumann control, Mixed control...), \iain{see equation \eqref{eq:dphs-balance}}. \iain{In fact, boundary conditions are required to describe a dPHS even if not explicitly appearing in equation \eqref{eq:dphs-balance}.} \textcolor{black}{An intrinsic difficulty associated with dPHS is their high dimensionality when discretized, which requires special care even in one-dimensional dPHS such as the one illustrated in Section \ref{sec:num}.}

\subsection{GP-dPHS} \label{section:Gp-dpHs}

We can now introduce GP-dPHS. Consider a dPHS whose Hamiltonian functional $\mathcal{H}$ is partially or totally unknown. We want to learn an approximation of the Hamiltonian $\widehat{\mathcal{H}}$. For this, we \iain{follow the} GP-pHs method \citep{Beckers} by \iain{setting} a function-indexed GP prior $\widehat{\mathcal{H}}$ on the \iain{true} Hamiltonian functional $\mathcal{H}$, yielding $\widehat{\mathcal{H}}\sim\GP(m(\balpha),k(\balpha,\balpha'))$, where $m : \mathcal{X}\rightarrow \R$ and $k : \mathcal{X}\times\mathcal{X}\rightarrow \R$ are prescribed mean and kernel functions. Let's also suppose that operators $\mathcal{J}$ and $\mathcal{R}$ are also unknown. We replace them by parametrized operators $\widehat{\mathcal{J}}$ formally skew-symmetric and $\widehat{\mathcal{R}}$ formally symmetric. We note $\widehat{\mathcal{J}_R} := \widehat{\mathcal{J}}-\widehat{\mathcal{R}}$. The resulting approximate dPHS is
\begin{equation}\label{eq:approx-dphs}
        %\dot{\balpha}(\balpha) 
        %f_{\dot{\balpha}}(\balpha) 
        \dot{\balpha} = \widehat{\mathcal{J}_R}\,\delta_{\balpha} \widehat{\mathcal{H}} (\balpha)  \ \ \ + \ \ \ \text{boundary conditions.}
\end{equation}
%coupled with the boundary condition that can also be supposed unknown and modelled. 
We denote $f_{\dot{\balpha}}(\balpha) :=\widehat{\mathcal{J}_R}\,\delta_{\balpha} \widehat{\mathcal{H}} (\balpha)$, \iain{the right-hand side of \eqref{eq:approx-dphs} viewed as a function of $\balpha$. From the stability of GPs via the linear map $\widehat{\mathcal{H}} \mapsto \widehat{\mathcal{J}_R}\,\delta_{\balpha} \widehat{\mathcal{H}}$, $f_{\dot{\balpha}}(\balpha)$ is an $\mathcal{X}$-valued GP (Section \ref{sec:gp}) with}
\begin{align*}
    f_{\dot{\balpha}}(\balpha) \sim \GP(\widehat{\mathcal{J}_R}\,\delta_{\balpha} m(\balpha), \widehat{\mathcal{J}_R}\,(\delta_{\balpha,\balpha'}k(\balpha,\balpha'))\,\widehat{\mathcal{J}_R^*}).
\end{align*} 
%By stability via the linear map $\widehat{\mathcal{H}} \mapsto \widehat{\mathcal{J}_R}\,\delta_{\balpha} \widehat{\mathcal{H}}$, $f_{\dot{\balpha}}(\balpha)= \widehat{\mathcal{J}_R}\,\delta_{\balpha} \widehat{\mathcal{H}} (\balpha)$ is a function-valued GP
Then, given $m_{\dot{\balpha}}$, the posterior expectation of the GP $ f_{\dot{\balpha}}(\balpha)$ conditioned on measurement data, %and given a sampled vector field $\widetilde{f}_{\dot{\balpha}}(\balpha,\omega),\ \omega \in \Omega $ from the posterior distribution of $f_{\dot{\balpha}}$, 
we can solve the differential equation (coupled with boundary conditions) $\dot{\balpha}(t) = m_{\dot{\balpha}}(\balpha(t))$ with $\balpha(0) = \balpha_0 \in \mathcal{X}.$ Alternatively we can replace $m_{\dot{\balpha}}$ with a sample of the posterior GP. \iain{Note that unlike GP-dPHS as described in \citep{tan24}, this version of GP-dPHS is still infinite-dimensional at this stage. Its discretization is achieved using the Partitioned Finite Element Method (PFEM) in Section \ref{sec:prior-pfem-gp-dphs}. This procedure consisting in first studying/controlling a dPHS and then discretizing is known as late lumping, as opposed to first discretizing it and then controlling it (early lumping), see \citep{marko2018early}.} \textcolor{black}{Finally, this approach ensures \emph{discretization invariance} over the estimation of $\mathcal{H}$ and the related dynamics, as described in e.g. \citep{lassas, stuart_bayes_inverse_pb}.}

\section{Discretization of the GP-dPHS with PFEM : PFEM-GP-dPHS}\label{sec:discretization-1D}

To perform numerical simulations, we need to discretize the GP $f_{\dot{\balpha}}$, or more precisely, write a GP $f_{\underline{\dot{\alpha}}}$ on the discrete coefficients of the discretization of $\dot{\balpha}$. To achieve this, we combine the GP-dPHS method with the partitioned finite element method (PFEM) \citep{PFEM} which has the advantage of being \denis{structure-preserving, in so far as the power balance of the dPHS is preserved at the discrete level}, see equation \eqref{eq:dphs-balance-wave}. We call this the PFEM-GP-dPHS method, which we illustrate below on a dissipative nonlinear wave equation in 1D.

\subsection{Dissipative nonlinear wave equation in dPHS form}

The first step is to write the 1D wave equation in dPHS form. This equation with Neumann boundary control is given by ($\O = [0,\ell]$) for $t \ge 0, x \in \O$

\begin{equation}
    \left\lbrace
    \begin{array}{rcl}
    \rho(x)\,\partial_{tt}^2 w(t,x) - \partial_x \left( s(x,\partial_x w(t,x))\,\partial_x w(t,x) \right) + \nu(x)\, \partial_{t} w(t,x)  &=& 0, \\
     -s(0,\partial_x w(t,0)) \partial_x w(t,0)  &=&  u_L(t), \\
    s(\ell,\partial_x w(t,\ell)) \partial_x  w(t,\ell) &=&  u_R(t).
    \end{array}
    \right.
\end{equation}

\noindent
We set the energy variable $\alpha_q := \partial_x w $ the strain and $\alpha_p := \rho \partial_t w $ the linear momentum. The Hamiltonian $\mathcal{H}$ of this system can be defined by  $\mathcal{H}(\alpha_q, \alpha_p) := \frac{1}{2} \int_0^\ell \int_0^{\alpha_q(t,x)} s(x,\alpha_q')\alpha_q'\,d\alpha_q'  + \frac{\alpha_p(t,x)^2}{\rho(x)} {\rmd}x.$ Computing the variational derivatives of $\mathcal{H}$ with respect to the energy variables leads to the co-energy variables $e_q := \delta_{\alpha_q} \mathcal{H} = s(.,\alpha_q)\alpha_q $, the stress, and $e_p := \delta_{\alpha_p} \mathcal{H} = \rho^{-1}\alpha_p$, the velocity. With Newton's second law and Schwarz's lemma, we get  the port-Hamiltonian system representing a damped vibrating string with Neumann boundary control and Dirichlet boundary observation:
\begin{equation}\label{eq:the_dpHs}
    \begin{pmatrix} \partial_t \alpha_q \\ \partial_t \alpha_p \end{pmatrix}
    =
    \underbrace{\begin{bmatrix} 0 & \partial_x \\ \partial_x & -\nu\end{bmatrix}}_{\mathcal{J_R}=\mathcal{J}-\mathcal{R}}
    \begin{pmatrix} e_q \\ e_p \end{pmatrix},\qquad \left\lbrace
    \begin{array}{rcl}
     e_q(t,0) &=& u_L(t), \ \ \ \ \  y_L(t) = e_p(t,0), \\
    e_q(t,\ell) &=& u_R(t), \ \ \ \ \  y_R(t) = e_p(t,\ell).
    \end{array}
    \right.
\end{equation}
\noindent
The following power balance can be computed for this system:
\begin{equation}\label{eq:dphs-balance-wave}
\frac{d}{dt} \mathcal{H}(t)= -\int_0^\ell \nu(x)\,{e_p}(t,x)^2   {\rmd}x + u_L(t) y_L(t) + u_R(t) y_R(t)\,.
\end{equation}

\subsection{The Partitioned Finite Element Method (PFEM)}\label{sec:pfem}

Let us apply the partitioned finite elements method (PFEM) \citep{PFEM} on the above example. Let $\varphi_q$ and $\varphi_p$ be smooth test functions. %, and $\delta_{mx}$ denote the Kronecker symbol. 
One can write the weak formulation and integrate the second line by parts to make the controls $u_L$ and $u_R$ appear:
\begin{equation*}
    \left\lbrace
    \begin{array}{rcl}
    \int_0^\ell  \dot{ \alpha_q}(t,x) \varphi_q(x) {\rmd}x &=& \int_0^\ell \partial_x e_p(t,x) \varphi_q(x) {\rmd}x, \\
    \int_0^\ell  \dot{ \alpha_p}(t,\cdot) \varphi_p {\rmd}x &=& - \int_0^\ell e_q(t,.) \partial_x \varphi_p {\rmd}x -\int_0^\ell \nu e_p(t,.)  \varphi_p {\rmd}x + u_R(t) \varphi_p(\ell) + u_L(t) \varphi_p(0), \\
    y_L(t) &=& e_p(t,0), \\
    y_R(t) &=& e_p(t,\ell).
    \end{array}
    \right.
    \label{eq:PFEM_pur}
\end{equation*}

Now, let $(\varphi_q^i)_{1 \le i \le N_q}$ and $(\varphi_p^j)_{1 \le j \le N_p}$ be two finite element families of approximation functions  for $\alpha_q$ and $\alpha_p$ respectively, so that $ \alpha_q(t,x) \simeq  \underline{\varphi_q }^\top(x) \underline{\alpha_q}(t)$, $\alpha_p(t,x) \simeq \underline{\varphi_p }^\top(x) \underline{\alpha_p}(t)$ (likewise for $\dot{\alpha_q}$ and $\dot{\alpha_p}$) and $ e_q(t,x) \simeq  \underline{\varphi_q }^\top(x) \underline{e_q}(t)$, $e_p(t,x) \simeq \underline{\varphi_p }^\top(x) \underline{e_p}(t)$. We next replace the variables by their approximations, and write the system in matrix form
% Suipprimer gL gR
\begin{equation}
    \left\lbrace
    \begin{array}{rcl}
    \underbrace{\begin{bmatrix}
    M_q & 0  \\
    0 & M_p 
    \end{bmatrix}}_{M} 
     \underbrace{\begin{pmatrix}
    \underline{\dot{\alpha}_q}(t) \\
    \underline{\dot{\alpha}_p}(t)
    \end{pmatrix}}_{ \underline{\dot{\alpha}}}
    &=&
   \underbrace{ \begin{bmatrix}
    0 & D &  \\
    -D^\top & -R_{22} 
    \end{bmatrix}}_{J-R}
     \underbrace{ \begin{pmatrix}
    \underline{e_q}(t) \\
    \underline{e_p}(t) 
    \end{pmatrix}}_{ \underline{e}} + \underbrace{ \begin{bmatrix}
    0 & 0 \\
    B_L & B_R 
    \end{bmatrix}}_{G}
   \underbrace{  \begin{pmatrix}
   u_L(t) \\
    u_R(t) 
    \end{pmatrix}}_{\underline{u}}, \\
   \underline{y}(t) \coloneqq  {\begin{pmatrix}
     y_L(t) \\
     y_R(t)
    \end{pmatrix}}
    &=&
   \begin{bmatrix}
   
    0 & B_L^\top \\
    0 & B_R^\top 
    \end{bmatrix}
    \begin{pmatrix}
    \underline{e_q}(t) \\
    \underline{e_p}(t) 
    \end{pmatrix},
    \end{array}
    \right.
    \label{eq:PFEM}
\end{equation}
with square matrices $ (M_q)_{ij} := \int_0^\ell \varphi_q^j(x) \varphi_q^i(x) {\rmd}x$,~ $ (M_p)_{k\ell} := \int_0^\ell \varphi_p^\ell(x) \varphi_p^k(x) {\rmd}x$,~ $ (R_{22})_{k\ell} := \int_0^\ell \nu(x) \varphi_p^\ell(x) \varphi_p^k(x) {\rmd}x$, and rectangular matrices~ $(D)_{i\ell} := \int_0^\ell \partial_x \varphi_p^\ell(x) \varphi_q^i(x) {\rmd}x$,~ $(B_L)_{k} := \varphi_p^k(0)$ and ~$ (B_R)_{k} := \varphi_p^k(\ell)$.
To close this system, we need to approximate the constitutive relations. Write them in their weak form:

\begin{equation}
    \left\lbrace
    \begin{array}{rcl}
    \int_0^\ell e_q(t,x) \varphi_q(x) {\rmd}x &=& \int_0^\ell s(x,\alpha_q(t,x)) \alpha_q(t,x) \varphi_q(x) {\rmd}x, \\
    \int_0^\ell e_p(t,x) \varphi_p(x) {\rmd}x &=&  \int_0^\ell \frac{1}{\rho(x)} \alpha_p(t,x) \varphi_p(x) {\rmd}x.
    \end{array}
    \right.
    \label{eq:constitutive_relation}
\end{equation}

\noindent
Then the matrix form of the discrete weak formulation of the constitutive relations is

\begin{equation}
       \left\lbrace
    \begin{array}{rcl}
    M_q \underline{e_q}(t) &=& M_s \underline{\alpha_q}(t), \\
    M_p \underline{e_p}(t) &=& M_\rho \underline{\alpha_p}(t),
    \end{array}
    \right.
    \label{eq:matrix_constitutive_relation}
\end{equation}
 
\noindent
where $(M_s)_{ij} := \int_0^\ell s(x,\alpha_q(t,x)) \varphi_q^j(x) \varphi_q^i(x) {\rmd}x$ ~and~ $(M_\rho)_{k\ell} := \int_0^\ell \frac{1}{\rho(x)} \varphi_p^\ell(x) \varphi_p^k(x) {\rmd}x$.

To retrieve a finite-dimensional pHs in equation~\eqref{eq:PFEM}, we can define the discrete Hamiltonian $\mathcal{H}^d$  as the evaluation of $\mathcal{H}$ on the approximation of the state variables
$
    \mathcal{H}^d(\underline{\alpha_q}, \underline{\alpha_p}) := \frac{1}{2} \underline{\alpha_q}^\top M_s(\alpha_q) \underline{\alpha_q} + \frac{1}{2} \underline{\alpha_p}^\top M_\rho \underline{\alpha_p}
$, and can write that

\begin{equation}
    \begin{bmatrix} \underline{e_q} \\ \underline{e_p} \end{bmatrix} =M^{-1}  \begin{bmatrix} \nabla_{\underline{\alpha_q}}   \\ \nabla_{\underline{\alpha_p}} \end{bmatrix} \mathcal{H}^d(\underline{\alpha_q},\underline{\alpha_p}) .
\end{equation}
\iain{The main perk of PFEM is that, whatever the chosen mesh refinement, we have the exact discrete power balance given by}
\denis{(compare with equation \eqref{eq:dphs-balance-wave}):
$$\frac{d}{dt} \mathcal{H}^d(t )= - \underline{e_p}^\top(t) R_{22}\,\underline{e_p}(t) + u_L(t) y_L(t) + u_R(t) y_R(t)\,.
$$
Numerical simulations of many dPHS models using PFEM can be performed using the \texttt{SCRIMP} environment \citep{FERRARO2024119}, see {\texttt{https://g-haine.github.io/scrimp/} }.
} \textcolor{black}{In Appendices \ref{subsec:2Ddphs} and \ref{subsec:2DPFEM}, we present an example of a 2D dPHS along with its discretization using~PFEM.}

\subsection{Expressing PFEM as an affine transformation of the Hamiltonian functional}\label{sec:pfem_affine}

\iain{Assume} that the Hamiltonian functional $\mathcal{H}$ \iain{is unknown}, and that the parameter $\nu$ \sylvain{is} \iain{given}. A first way to \sylvain{identify} \iain{$\mathcal{H}$} would be to use the GP-PHS method from \citep{Beckers} by \iain{setting} a GP prior on %\iain{a previously chosen} 
\iain{some} discretized Hamiltonian $\mathcal{H}^d$ (early lumping). \iain{Unfortunately}, this approach \sylvain{leads to} a number of hyperparameters proportional to the number of degrees of freedom used for the discretization as we would have to define a set of hyperparameter for each degree of freedom. 
\sylvain{
To solve this issue, we \iain{use the} GP-dPHS method by \iain{first setting} a function-indexed GP prior $\widehat{\mathcal{H}}$ on the \iain{true} Hamiltonian functional $\mathcal{H}$ itself~: $\widehat{\mathcal{H}}\sim\GP(m(\balpha),k(\balpha,\balpha'))$ (late lumping), see Section \ref{section:Gp-dpHs}.} \iain{ $f_{\dot{\balpha}}(\balpha)= \widehat{\mathcal{J}_R}\,\delta_{\balpha} \widehat{\mathcal{H}} (\balpha)$ is a function-valued GP, see Section \ref{section:Gp-dpHs}. Explicit mean and kernel are described in Section \ref{sec:prior-pfem-gp-dphs}.} %\iain{In section \ref{sec:prior-pfem-gp-dphs}, we describe the resulting GP on $f_{\dot{\balpha}}(\balpha)= \widehat{\mathcal{J}_R}\,\delta_{\balpha} \widehat{\mathcal{H}} (\balpha) $}.
%To keep some freedom on the way we define parameters, we want to define them as function of the space variable. 
\iain{We next choose the hyperparameters of $m$ and $k$ as continuous functions themselves, which enables us to discretize them independently from the PFEM discretization of a given dPHS, see Section \ref{sec:prior-pfem-gp-dphs}. We will use piecewise linear approximations but other choices can be considered such as splines, trigonometric sums or wavelets.}
%For that \sylvain{ we follow} the GP-dpHs method by applying a GP prior to the Hamiltonian functional. Following the section \ref{section:Gp-dpHs}, we thus replace $\mathcal{H}$ by an approximate $\widehat{\mathcal{H}}\sim\GP(m(\balpha),k(\balpha,\balpha'))$, 
 %The unknown function $\nu$ is modelled by a parameterized function $\widehat{\nu}$ which leads to the parametrize operator $\widehat{\mathcal{J}_R} = \left[[0, \partial_x]^\top,[\partial_x,\widehat{\nu}]^\top\right]$. 

%\sylvain{Séparer la partie continue de la partie discrète par un nouveau paragraphe}
\iain{We then apply the PFEM method on the dPHS $\dot{\balpha} =f_{\dot{\balpha}}(\balpha)$. We show below that PFEM can be viewed as applying a linear transformation on $f_{\dot{\balpha}}(\balpha)$, thus yielding a GP $f_{\dot{\underline{\alpha}}}(\underline{\alpha},\underline{u})\in\R^{N_q+N_p}$ which is vector-valued but also vector-indexed by the \emph{pair} $(\underline{\alpha},\underline{u})\in \R^{N_q+N_p}\times \R^2$; see Section \ref{sec:pfem} for the introduction of $\underline{\alpha}$ and $\underline{u}$.}
%We then apply the PFEM method to the approximate dpHs $\dot{\balpha} =\mathcal{J}_R\bm{e} (\balpha)$, with $\bm{e} :=\delta_{\balpha}\widehat{\mathcal{H}}$, which for this case give back equation~\eqref{eq:PFEM}.\\
%Then, to have a discretized GP, we want to have a GP on $\underline{\dot{\alpha}}$, and for that, 

We thus express $\underline{\dot{\alpha}}$ as combination of affine transformations applied to $\widehat{\mathcal{H}}$. For this, we start from the left-hand side of equations~\eqref{eq:constitutive_relation} and~\eqref{eq:matrix_constitutive_relation} to link $\underline{e}$ and $\widehat{\mathcal{H}}$, yielding

\begin{equation}
    M\underline{e}(\balpha) = \int_0^\ell\Phi(x) (\delta_{\balpha} \widehat{\mathcal{H}})(\balpha)(x)  {\rmd}x\,, \ \ \ \Phi(x) = \begin{bmatrix} \underline{\varphi_q}(x) &  0_{N_q,1}\\ 0_{N_p,1} & \underline{\varphi_p}(x) \end{bmatrix} .
    \label{eq:transfo}
\end{equation}

\noindent Above we have $\Phi \in \R^{(N_q+N_p)\times2}$. We next substitute $\balpha$ with its finite dimensional approximation $\Phi^T\underline{\alpha}\simeq \balpha$ and use equation~\eqref{eq:PFEM} which links $\dot{\underline{\alpha}}$ and $\underline{e}$, to find 
\begin{equation}
    f_{\dot{\underline{\alpha}}}(\underline{\alpha},\underline{u}) = \dot{\underline{\alpha}} =M^{-1}(J-R)M^{-1}  \int_0^\ell\Phi(x) (\delta_{\balpha} \widehat{\mathcal{H}}) (\Phi^T\underline{\alpha})(x) {\rmd}x  +G\underline{u}\,.
    \label{eq:affine_transfo}
\end{equation}

\noindent Equation \eqref{eq:affine_transfo} thus expresses $f_{\dot{\underline{\alpha}}}(\underline{\alpha},\underline{u})$ as an affine transformation of $\widehat{\mathcal{H}}$. and we can derive the GP followed by $ f_{\dot{\underline{\alpha}}}(\underline{\alpha},\underline{u})$. The prior of this so called PFEM-GP-dPHS depends on the hyperparameters of the prior $\widehat{\mathcal{H}}$, as well as the parameter $\nu$. As a function of $\underline{\alpha}$, the discrete output $\underline{y}$ is also a GP~since
\begin{equation}
   \underline{y}(\underline{\alpha}) =G^\top M^{-1} \int_0^\ell\Phi(x) (\delta_{\balpha} \widehat{\mathcal{H}}) (\Phi^T\underline{\alpha})(x) {\rmd}x \,.
\end{equation}

After the GP $\widehat{\mathcal{H}}$ is conditioned on measurement data, we obtain a posterior vector-valued GP whose (posterior) expectation is denoted $m_{\dot{\underline{\alpha}}}(\underline{\alpha},\underline{u})$. Finally, we emulate the true dPHS by solving $\underline{\dot{\alpha}}(t) =  m_{\dot{\underline{\alpha}}}(\underline{\alpha}(t),\underline{u}(t))$  with $\underline{{\alpha}}(0) = \underline{{\alpha}}_0 \in \mathcal{X}$.% where $\underline{\alpha}\mapsto \widehat{f}_{\dot{\underline{\alpha}}}(\underline{\alpha},\underline{u},\omega)$ is a sample of the conditioned PFEM-GP-dpHs described above.

\section{Prior of the PFEM-GP-dPHS}\label{sec:prior-pfem-gp-dphs}

\subsection{Selecting a prior for the Hamiltonian functional}\label{subsec:select_prior_h}

%Now in order to obtain a prior for $ f_{\dot{\underline{\alpha}}}$, we need do define a prior on $\widehat{\mathcal{H}}$.
\sylvain{Following equation~\eqref{eq:affine_transfo} a prior on $\widehat{\mathcal{H}}$ should be defined to obtain the prior of $ f_{\dot{\underline{\alpha}}}$.}
For this \citep{NAPI-MPC} proposed using a zero mean and a squared exponential (SE) kernel. However, in most physical systems, the Hamiltonian usually goes to infinity when a variable goes to infinity. This behaviour is absolutely not captured by GP with zero mean and SE kernel, because sampled functions from a GP with a constant zero mean value and SE kernel are necessary bounded. To remedy this, we propose to use a GP prior with a quadratic mean and a SE kernel. This will first allow the sampled Hamiltonians to go to infinity when a variable goes to infinity and second, as in physical systems there are often quadratic dependencies of the Hamiltonian on some of the variables, the mean will approximate the quadratic part of the Hamiltonian, and the kernel will focus on the non-quadratic part of the Hamiltonian, allowing better modelling of the system. \iain{Note that the quadratic part of the Hamiltonian corresponds to the linear part of the dPHS, hence the kernel will mainly focus on approximating the nonlinear part of the dPHS.} Thus the scalar prior mean $m$ and kernel $k$ selected for the GP of $\widehat{\mathcal{H}}$ of our example, such that $\widehat{\mathcal{H}} \sim \GP(m,k)$, are given by

\begin{align} 
m(\balpha) &=  \int_0^\ell m_{q}(x)\alpha_q^2(x)\,dx +  \int_0^\ell m_{p}(x)\alpha_p^2(x)\,dx, 
\label{eq:mean_GP}\\
    k(\balpha,\balpha') &= \sigma_f^2 \exp\left[-\frac{1}{2}\left( \int_0^\ell \frac{(\alpha_q(x)-{\alpha_q'}(x))^2}{\ell_q^2(x)}\,dx +  \int_0^\ell \frac{(\alpha_p(x)-{\alpha_p'}(x))^2)}{\ell_p^2(x)}\,dx \right)\right], 
\label{eq:kernel_GP}
\end{align}

where the functions $m_{q}$, $m_{p}$, $\ell_q$ and $\ell_p$ are nonnegative hyperparameters of the GP $\widehat{\mathcal{H}}$, \sylvain{and $\sigma_f \in \mathbb{R}^+$ is also a hyperparameter of the GP.}  \iain{In Section \ref{sec:num} the functions $m_q, m_p, \ell_q$ and $\ell_p$ are discretized using piecewise linear functions; a discretization based on 3\textsuperscript{rd} order polynomials is studied in Appendix \ref{app:polyn}.} \textcolor{black}{Note that this quadratic mean is a form a physics-informed prior, as the Hamiltonians associated to many dPHS have a quadratic term, e.g. the one associated to the kinetic energy in classical and quantum mechanics. A more detailed study of the role and the impact of the quadratic mean is given in Appendix~\ref{app:role_quad_mean}.}

\subsection{Derivation of the prior over \texorpdfstring{$ f_{\dot{\underline{\alpha}}}$}{f-alpha-dot}  }

 To obtain the prior of $f_{\dot{\underline{\alpha}}}$, each affine transformation appearing in equation~\eqref{eq:affine_transfo} is applied to the GP of $\widehat{\mathcal{H}}$. First we apply the variational derivative which we decompose as $\delta_{\balpha} = [\delta_{\alpha_q} \ \delta_{\alpha_p} ]^\top$. This yields a GP prior over $\bm{e}(\balpha) \sim \GP(m_e(\balpha), k_e(\balpha,\balpha'))$ with $ m_e(\balpha) :=  \delta_{\balpha} m(\balpha)$ and $  k_e(\balpha,\balpha') :=   \delta_{\balpha} \delta_{\balpha'}k(\balpha,\balpha')$. Explicitly we have the two formulas

\begin{equation}
    m_e(\balpha) = \begin{bmatrix}m_{q}\alpha_q \\ m_{p}\alpha_p  \end{bmatrix}\,,
\end{equation}

\begin{equation}
    \frac{k_e(\balpha,\balpha')}{ k(\balpha,\balpha')}  =   \begin{bmatrix}
   \frac{1}{\ell_q^2}\mathcal{I} -\frac{1}{\ell_q^2}(\alpha_q -\alpha_q') \otimes \frac{1}{\ell_q^2}(\alpha_q-\alpha_q')  & -\frac{1}{\ell_q^2}(\alpha_q -\alpha_q') \otimes \frac{1}{\ell_p^2}(\alpha_p-\alpha_p') \\-
    \frac{1}{\ell_p^2}(\alpha_p -\alpha_p') \otimes \frac{1}{\ell_q^2}(\alpha_q-\alpha_q')& \frac{1}{\ell_p^2}\mathcal{I} - \frac{1}{\ell_p^2}(\alpha_p -\alpha_p') \otimes \frac{1}{\ell_p^2}(\alpha_p-\alpha_p') 
    \label{eq:consti_kernel}
\end{bmatrix}\,,
\end{equation}

Note that above, we have a function-valued mean and a symmetric p.s.d. operator-valued kernel. \iain{Following Section \ref{sec:pfem_affine}}, multiply by $\Phi$, integrate and replace $\balpha$ by its finite dimensional approximation $\Phi^\top\underline{\alpha}$. This yields a GP prior over $\underline{e}(\underline{\alpha})$ given by $\underline{e}(\underline{\alpha}) \sim \GP(m_{\underline{e}},k_{\underline{e}})$ where $m_{\underline{e}}(\underline{\alpha}) = M^{-1}\int_0^\ell  \Phi(x)   m_e(\underline{\alpha})(x) {\rmd}x$ and where the symmetric p.s.d matrix-valued kernel is 
\begin{align*}
    k_{\underline{e}}(\underline{\alpha},\underline{\alpha}') = M^{-1}\bigg(\int_0^\ell  \Phi(x)  ( k_e(\underline{\alpha},\underline{\alpha}') \Phi^\top)(x) {\rmd} x\bigg)M^{-1}\,.
\end{align*} 
Their expressions are given by
\begin{equation}\label{eq:mean_e}
   m_{\underline{e}}(\underline{\alpha}) = \begin{bmatrix} (M_{m,q} \underline{\alpha_q})^\top &  (M_{m,p} \underline{\alpha_p})^\top   \end{bmatrix}^\top,
\end{equation}

\noindent with $(M_{m,q})_{ij} :=  \int_0^{\ell}  m_{q}\varphi_q^i \varphi_q^j {\rmd}x$ and $(M_{m,p})_{ij} :=  \int_0^{\ell} m_{p} \varphi_p^i \varphi_p^j{\rmd}x$ and 

\begin{equation}\label{eq:kernel_e}
        \frac{k_{\underline{e}}(\underline{\alpha}, \underline{\alpha}')}{k(\Phi^\top\underline{\alpha},\Phi^\top\underline{\alpha}')} =  M_k\begin{bmatrix}
       I_{N_q}-(\underline{\alpha_q}-\underline{\alpha_q'})(\underline{\alpha_q}-\underline{\alpha_q'})^\top  &-(\underline{\alpha_q}-\underline{\alpha_q'})(\underline{\alpha_p}-\underline{\alpha_p'})^\top  \\-(\underline{\alpha_p}-\underline{\alpha_p'})(\underline{\alpha_q}-\underline{\alpha_q'})^\top &I_{N_p}-(\underline{\alpha_p}-\underline{\alpha_p'})(\underline{\alpha_p}-\underline{\alpha_p'})^\top
       \end{bmatrix}M_k^\top,
\end{equation}

\noindent with $(M_{k,q})_{ij} :=  \int_0^\ell  \ell_q^{-2}\varphi_q^i \varphi_q^j{\rmd}x$, $(M_{k,p})_{ij} :=  \int_0^\ell \ell_p^{-2} \varphi_p^i \varphi_p^j{\rmd}x$ and $M_k = M^{-1}\text{Diag}(M_{k,q},M_{k,p})$. Then from here, to obtain the GP prior over $ f_{\dot{\underline{\alpha}}}$, we  apply the matrix multiplication $(J-R)$ to the GP $\underline{e}$ and add the mean $G\underline{u}$ so that finally
\begin{equation}
    f_{\dot{\underline{\alpha}}}(\underline{\alpha},\underline{u}) \sim \GP(M^{-1}(J-R) m_{\underline{e}}(\underline{\alpha}) +G\underline{u}\,,\ M^{-1}(J-R) k_{\underline{e}}(\underline{\alpha}, \underline{\alpha}') (J-R)^\top M^{-1}).
\end{equation}

\textcolor{black}{The results given in Sections \ref{sec:discretization-1D}, \ref{sec:prior-pfem-gp-dphs} and \ref{sec:num} are generalized to a 2D dPHS in Appendix \ref{app:2d}.} %\ref{subsec:PFEM-GP_2D}, \ref{subsec:PFEM_GP_nD}
%\iain{dire qu'on étend l'approche de Beckers à la dim infinie.}
%\iain{puis, nouveauté : discrétisation PFEM (utiliser le savoir faire EDP)}
%\iain{un des intérêts : on a un choix + flexible sur les fonctions de base pour représenter les paramètres}

%\paragraph{Efficient discretisation of the parameters : hyperparameters + physical parameters}

%Now that we have our finite-dimensional PFEM-GP-dpHs, we have also the parameters of $m$ and $k$ also $\nu$ that are function of the space variable. From here, anything is possible to discretize. One choice would be to decompose each function on a new basis of FE, defined on a bigger mesh to reduce the number of parameters compared to the classical GP-pHs method.

\section{Numerical results}\label{sec:num}

%\textcolor{black}{Examples of numerical results of PFEM-GP-dPHS for a 2D system are given in Appendix \ref{subsec:2Dnum}.}
We test PFEM-GP-dPHS on the dPHS described by equation~\eqref{eq:the_dpHs}. We decompose the Young's Modulus term $s$ into the sum of a quadratic term $T(x) = 2 - 4x(1-x)$ and a spatially modulated non-quadratic term $c(x)f(\alpha_q(t,x))$. Here, $c(x) = 2x(x-1)^2$ is the spatial modulator and $f(u) = \exp(-u^2)$ is the non-quadratic function. We thus have $s(x,\alpha_q(t,x)) = T(x) + c(x)f(\alpha_q(t,x))$. The linear density is given by $\rho(x) = 3 -  2.5x^2$. The Hamiltonian functional can be written as

\begin{equation}\label{eq:hamiltonian_simu}
   \mathcal{H}(\balpha) = \mathcal{H}(\alpha_q, \alpha_p) := \frac{1}{2} \int_0^\ell c(x)\int_0^{\alpha_q(t,x)} f(\alpha_q')\alpha_q'\,d\alpha_q' + T(x) \alpha_q^2(t,x)  + \frac{\alpha_p(t,x)^2}{\rho(x)}  {\rmd}x.
\end{equation}

Above, we have two quadratic terms depending on the parameters $T$ and $\rho^{-1}$ and one non-quadratic term depending on the parameter $c$. For simplicity, we take $\nu = 0$ so that only the Hamiltonian parameters $T$, $\rho$ and $c$ are supposed unknown. In the initial state, we set $\underline{\alpha}_0 = \underline{0}$, and we use $\underline{u}(t) = (u_L(t) ,u_R(t)) = (\sin(\pi t), 0)$.
$\alpha_q$ and $\alpha_p$ are discretized on a 21 or 41-point uniform mesh on $\mathbb{P}^1$ Lagrange elements, resulting in 42 or 82 degrees of freedom ($N_q = N_p \in \{21,41\}$), \textcolor{black}{resulting in a finite-dimensional PHS with high dimension $d\in\{42,82\}$}. The hyperparameters $m_q$, $m_p$, $\ell_q^{-1} $ and $\ell_p^{-1}$ are discretized on a mesh of $\mathbb{P}^1$ Lagrange elements with mesh size of $\Delta x \in \{0.1, 0.2\}$. Accounting for $\sigma_f$ and $\sigma_{\text{noise}}$, this yields a hyperparameter $\theta$ of size 26 or 46 (note that different meshes for each hyperparameter functions could be used). In Appendix \ref{app:polyn}, we approximate the hyperparameter functions with polynomials of degree 3, resulting in 18 hyperparameters. 
Note that using an early lumping approach would result in 166 hyperparameters. 

\iain{The training data is obtained by first simulating the dPHS described above for $t_f = 20s$ with the \texttt{SCRIMP} environment \citep{FERRARO2024119} with time step $\Delta t = 0.01s$ and number of simulation time steps $N_t = 2000$. In this single simulation} we select for $0 \leq t \leq 10 s$, $N_s = 35$ uniformly spaced time stamps $t_i $ ($0\leq i < N_s$) and the corresponding true $\underline{\alpha}_i$ and $\underline{\dot{\alpha}}_i$, resulting in $N_s\times (N_q+N_p)=1470$ training points. The estimation of the hyperparameters is done by minimizing the negative log marginal likelihood (NLML) using the L-BFGS-B method of \verb|scipy| along with the NLML gradient. %\iain{on a bounded domain}. 
After this minimization step, \iain{we estimate $f_{\dot{\underline{\alpha}}}$ with its conditional expectation with reference to the training data}. We next simulate the resulting dPHS using $m_{\underline{\dot{\alpha}}}$, the posterior mean of $ f_{\dot{\underline{\alpha}}}$ conditioned on the training data, by solving $\dot{\underline{\alpha}}(t) =m_{\underline{\dot{\alpha}}}(\underline{\alpha}(t),\underline{u}(t))$ with $\underline{\alpha}(0) = \underline{\alpha}_0$.
\textcolor{black}{A discussion on the computational cost of this procedure is given in Appendix \ref{subsec:computational_costs}.}

Given the Hamiltonian $\mathcal{H}(\balpha)$ in equation \eqref{eq:hamiltonian_simu} as well as the mean and kernel functions in equations \eqref{eq:mean_GP} and \eqref{eq:kernel_GP}, we expect the mean function $m(\balpha)$ to learn the quadratic part of $\mathcal{H}(\balpha)$, and the kernel $k(\balpha,\balpha')$ to learn the nonquadratic part of $\mathcal{H}(\balpha)$.

    \begin{figure}[t!]    
    \includegraphics[width=0.45\textwidth]{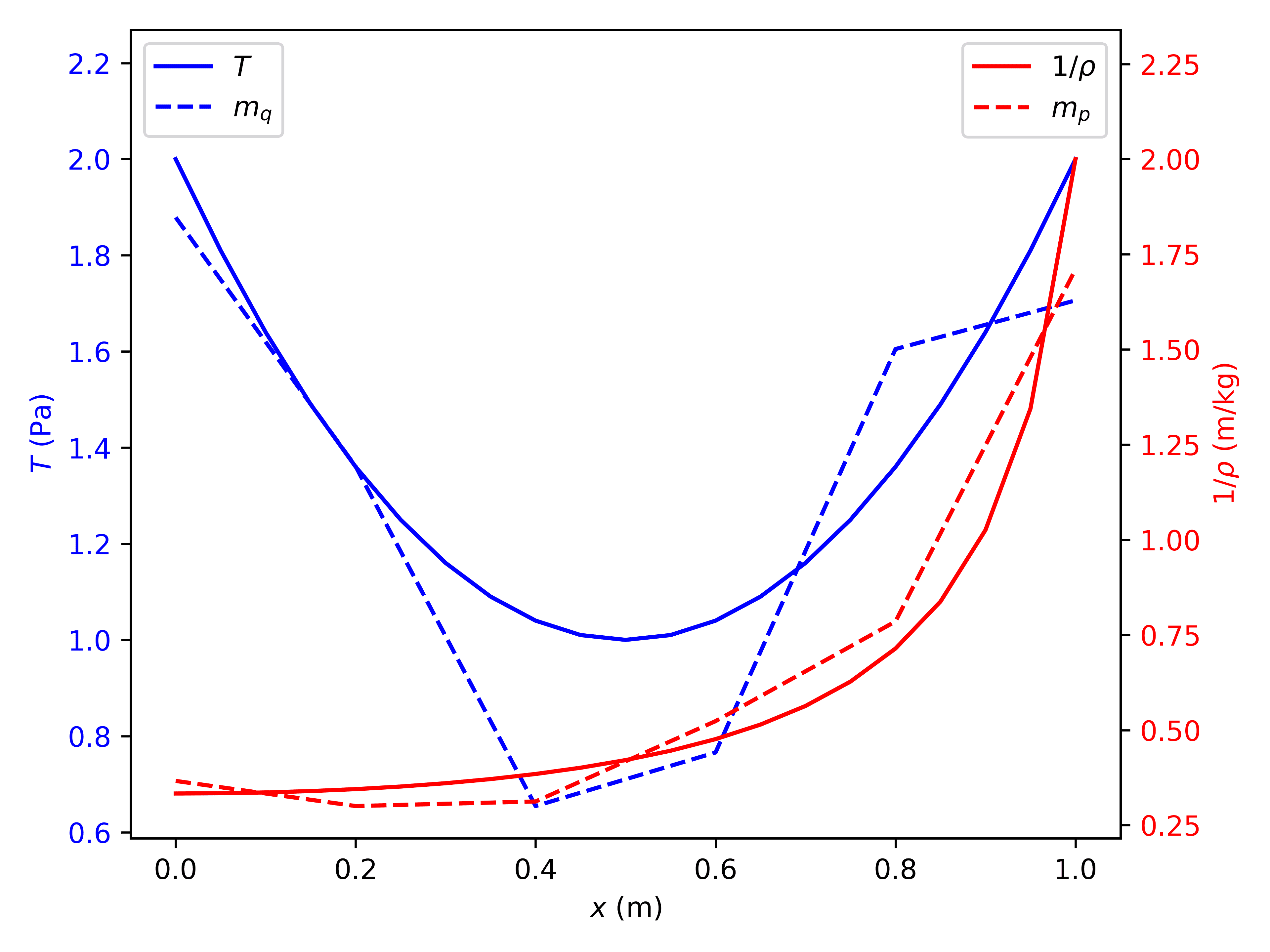}\centering
    \includegraphics[width=0.45\textwidth]{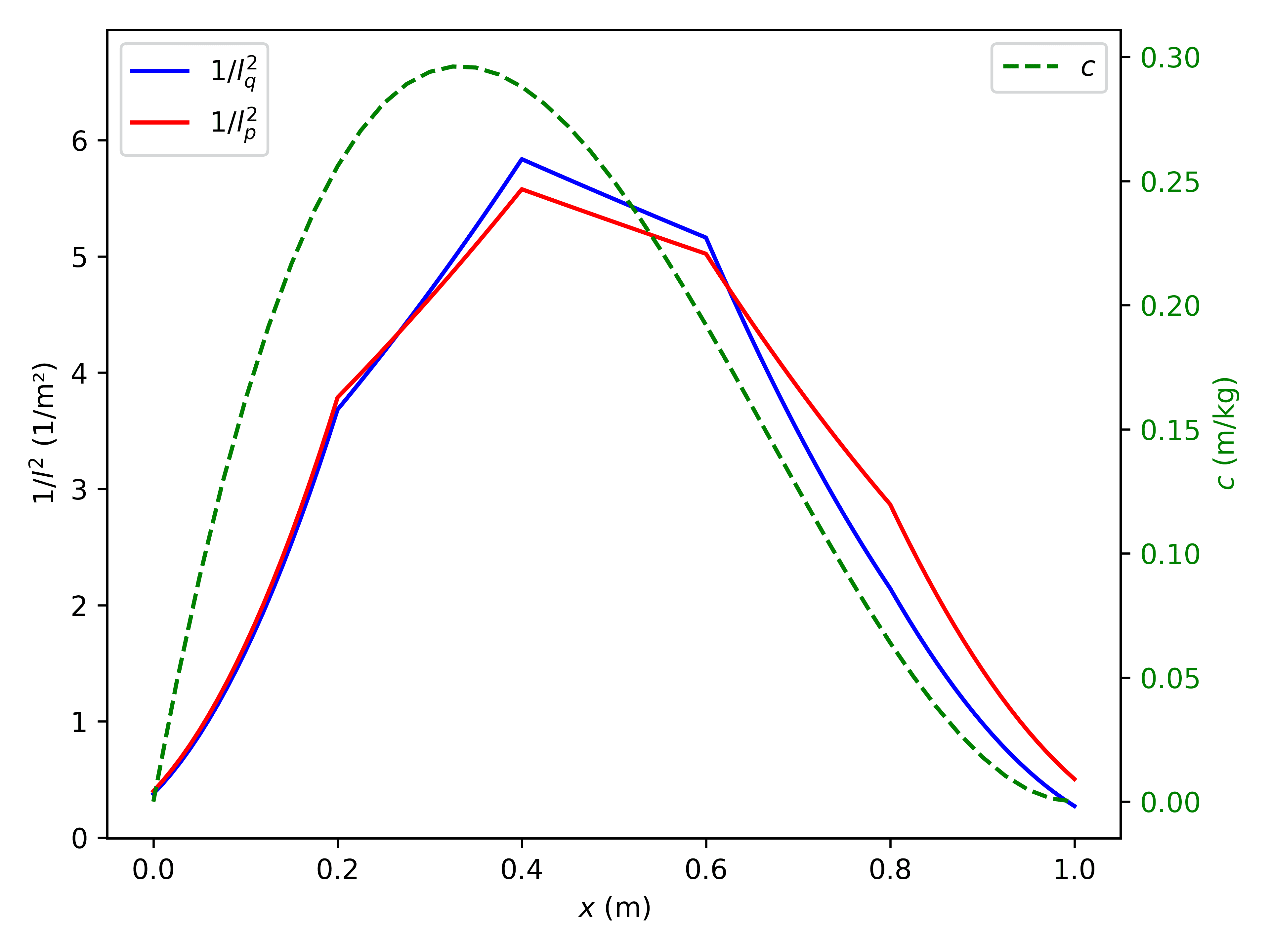}
    \caption{Left : estimation of the GP quadratic mean functions $m_q$ and $m_p$ using piecewise linear functions, versus true physical parameters $T$ and $\rho^{-1}$. Right : estimation of the GP functional lengthscales $\ell_q^{-2}$ and $\ell_p^{-2}$ using  piecewise linear functions, versus the function $c$. Here $N_q = N_p = 21$, $N_s = 35$ and the hyperparameter mesh size is $\Delta x = 0.2$ resulting in 26 hyperparameters.}
    \label{fig:hyperparameters}
\end{figure}

The estimated hyperparameter functions after successful NLML minimisation are given in Figure \ref{fig:hyperparameters}. Figure \ref{fig:hyperparameters} (left) shows that $m_q$ and $m_p$ approximately follow their respective counterparts $T$ and $\rho^{-1}$, indicating that the quadratic part of the true Hamiltonian has been approximated by the quadratic mean of the GP prior. Figure \ref{fig:hyperparameters} (right) shows  that the shape of the squared inverse correlation lengths $\ell_q^{-2}$ and $\ell_p^{-2}$ follows the shape of the spatial modulator $c$. We explain this result as follow. Interpreting $\ell_p(x)$ and $\ell_q(x)$ as lengthscales at ``dimension $x$'', we expect that the Hamiltonian is well approximated by the quadratic mean at locations (``dimensions $x$'') where its nonlinearity is mild, i.e. $c(x)\simeq 0$. Hence we expect that the tuned kernel $k$ should not contribute at such locations $x$. Ensuring this latter property corresponds to having infinite lengthscales $\ell_p(x)$ and $\ell_q(x)$, which is the case in Figure \ref{fig:hyperparameters} (right) at locations where $c(x)\simeq 0$. 

{Figure \ref{fig:alpha_p_mix} and {\ref{fig:alpha_q_mix} in Appendix \ref{subsec:additional_figs} describe simulated states using the approximate evolution equation $\dot{\underline{\alpha}}(t) =m_{\underline{\dot{\alpha}}}(\underline{\alpha}(t),\underline{u}(t))$, versus the \texttt{SCRIMP} output. Those figures provide visual evidence that on this example, the approximate dPHS $\dot{\underline{\alpha}}(t) =m_{\underline{\dot{\alpha}}}(\underline{\alpha}(t),\underline{u}(t))$ obtained using the PFEM-GP-dPHS framework provides a reasonable approximation of the true dPHS (equation \eqref{eq:dphs-balance})}.

Figure \ref{fig:rel_error_alpha_p} studies the behaviour of our method as a function of the discretization of the functions $m_q, m_p, \ell_q^{-1}$ and $\ell_p^{-1}$. First recall that we have at our disposal a single \texttt{SCRIMP} output $(\alpha_q^{\text{SCP}}, \alpha_p^{\text{SCP}})$ with fixed measurement points (up to 10$s$ in the simulation) as described above.~We~discretize the functions $m_q, m_p, \ell_q^{-1}$ and $\ell_p^{-1}$ on the same grid with a chosen mesh size $\Delta x \in \{0.01, 0.02, 0.025,$ $ 0.05, 0.1, 0.2, 0.25, 0.5\}$. Accounting for $\sigma_f$ and $\sigma_{\text{noise}}$, this yields a hyperparameter $\theta$ of size $2 + 4 N_{\Delta x}$, where $N_{\Delta x} = 1+\Delta x^{-1} \in \{101, 51, 41, 21, 11, 6, 5, 3\}$. Given this mesh size, we run $N_\text{opt} = 20$ L-BFGS-B NLML minimisations, with an initial guess chosen at random in $[1,2]^{2 + 4N_{\Delta x}}$. This yields 20 different optimized values $\theta_i^*$, $1\leq i \leq N_{\text{opt}}$. This in turn yields $N_{\text{opt}}$ posterior expectations $m^i_{\underline{\dot{\alpha}}}$, $N_{\text{opt}}$ approximated dPHS equations $\dot{\underline{\alpha}}(t) =m^i_{\underline{\dot{\alpha}}}(\underline{\alpha}(t),\underline{u}(t))$ and $N_{\text{opt}}$ solutions $(\alpha_q^i,\alpha_p^i)$. For each $i$, we compare $\alpha_p^i$ and the \texttt{SCRIMP} output $\alpha_p^{\text{SCP}}$ : at each time step $t_j \in \{0, \ldots, N_t\}$, we compute the relative $L^2$ spatial error $e_{i,j} = \|\alpha_p^i(t_j) - \alpha_p^{\text{SCP}}(t_j)\|_{L^2}/\|\alpha_p^{\text{SCP}}(t_j)\|_{L^2}$ (likewise for $\alpha_q$), and the spatial error averaged over time $e_i = \frac{1}{N_t}\sum_{j=0}^{N_t}e_{i,j}$ for $\alpha_p$. Figure \ref{fig:rel_error_alpha_p} (left) shows an example of $(e_{i,1}, \ldots, e_{i,N_t})$ for $\Delta x = 0.1$. Running this whole procedure we obtain, for each mesh size $\Delta x$, and a sample $(e_1, \ldots, e_{N_{\text{opt}}})$. Figure \ref{fig:rel_error_alpha_p} (right) depicts histograms of these data sets. This figure shows that our method works best when the hyperparameter mesh size is either $0.05$ or $0.1$, the former also being the mesh size of the reference \texttt{SCRIMP} simulation (21 mesh points on the spatial domain $[0,1]$). This figure also shows that the L-BFGS-B optimization is the most robust for those step sizes. For the other step sizes, an average relative error higher than $1$ is not uncommon, meaning that the null predictor performs better than PFEM-GP-dPHS in those cases. In Appendix \ref{app:posterior_variance}, we discuss the use of the posterior variance provided by PFEM-GP-dPHS for active learning purposes.
\begin{figure}[t!]
    \centering
     \includegraphics[width=0.45\linewidth]{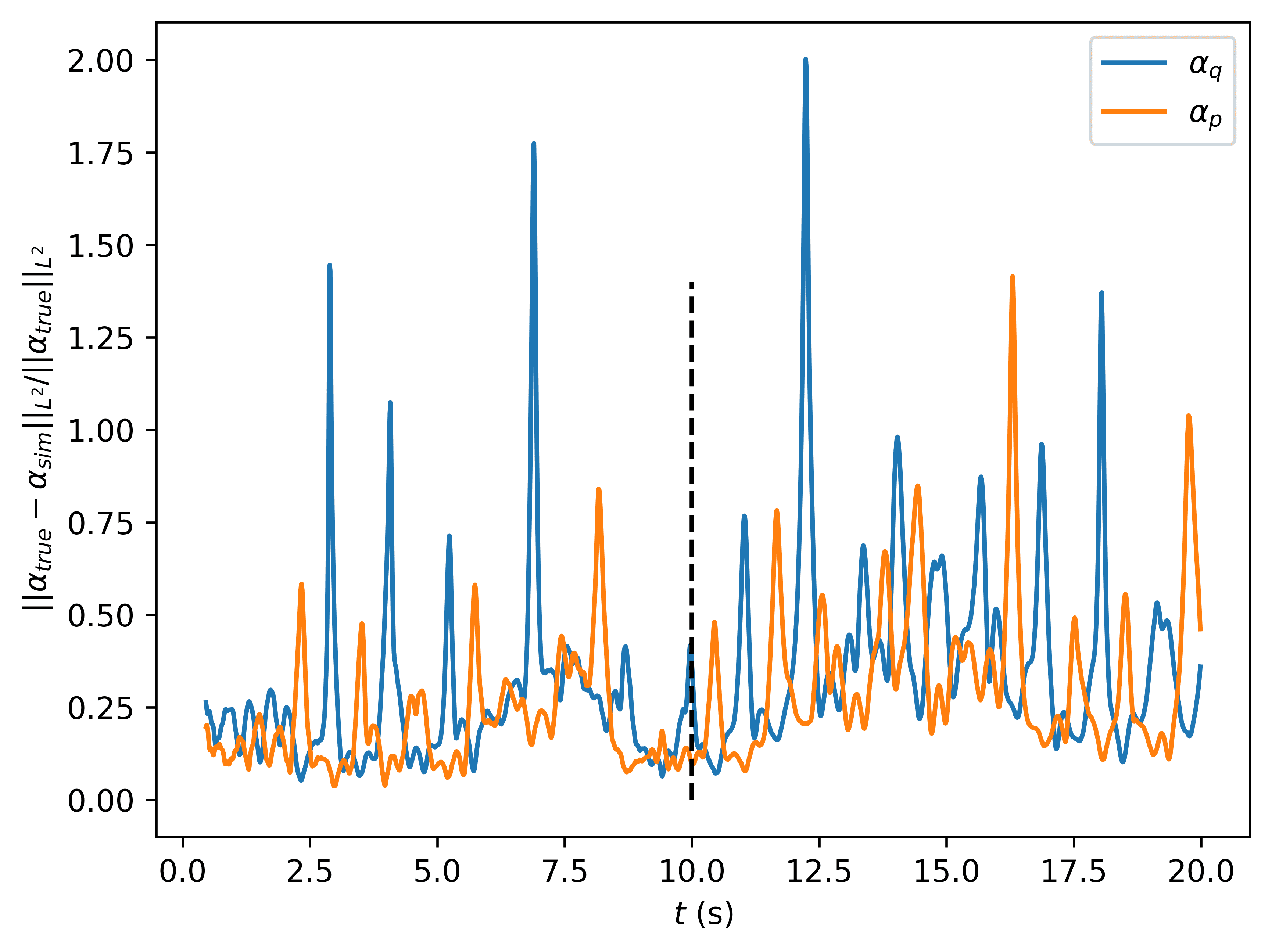}
    \includegraphics[width=0.45\linewidth]{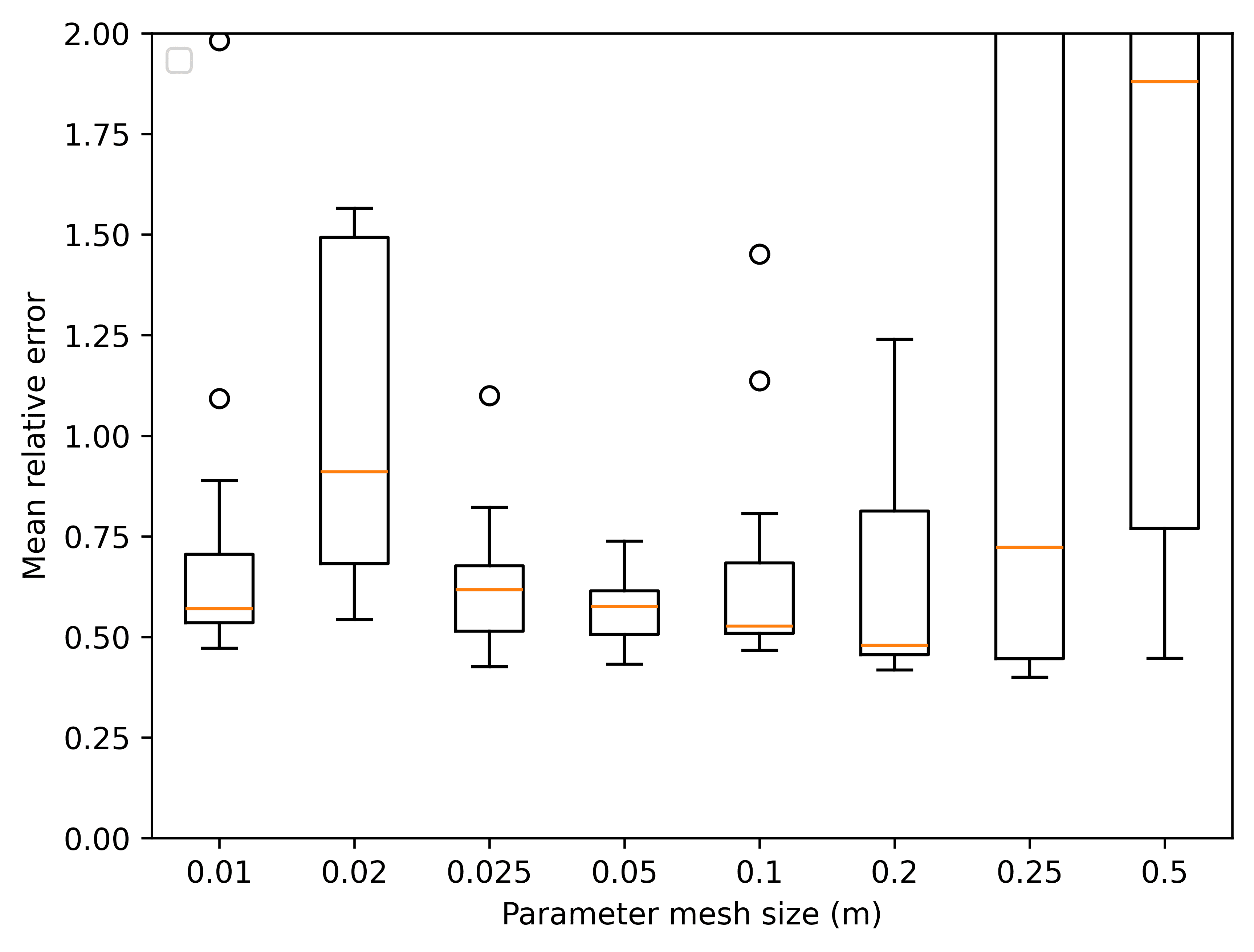}
    \caption{Study of the behaviour of PFEM-GP-dPHS as a function of the hyperparameter mesh size.}
    \label{fig:rel_error_alpha_p}
\end{figure}

\section{Conclusion}

By combining the physics-based prior of GP-dPHS and the PFEM discretization method, we were able to obtain a new GP prior that preserves the passivity of a given dPHS at the discrete level. Using a late lumping approach allowed us to introduce functional hyperparameters that can be later discretized independently from PFEM. By discretizing them as  piecewise linear functions, PFEM-GP-dPHS is shown to be capable of learning a 1D nonlinear wave equation. The choice of a quadratic mean in the GP prior over the Hamiltonian functional allowed us to better model the Hamiltonian, allowing the identification of quadratic components of the true Hamiltonian functional (linear part of the dPHS) and separating them from non-quadratic contributions.  Future works include  the further study of our framework in the presence of dissipation ($\nu\neq 0$) as well as nD dPHS, both cases being compatible with the \texttt{SCRIMP} Python package.  

%\textcolor{black}{Learning of 2D system using our method also work, and generally learning of dPHS of any dimensionality is theoretically possible.}
%\textcolor{black}{the use of the posterior variance for adaptative learning, and}

\subsection*{Acknowledgements} {Iain Henderson would like to thank the DRRP (Direction de la Recherche et des Ressources Pédago\-giques) of ISAE-Supaéro for funding the internship of Florian Courteville.}

\bibliographystyle{plainnat}

\bibliography{bibli.bib}

\clearpage

\appendix

\begin{center}
{\large 
\textbf{PFEM-GP-dPHS : a finite element framework for combining Gaussian processes and infinite-dimensional port-Hamiltonian systems \\ \vspace{0.5cm}
Supplemental material}}
\end{center}

\section{Study of the role of the quadratic mean}\label{app:role_quad_mean}
Kernel regression methods are universal approximators and thus the benefits of using a quadratic mean should be further examined.
\paragraph{Study of PFEM-GP-dPHS without using a quadratic mean}
We run the same workflow as the one we used to obtain Figure \ref{fig:rel_error_alpha_p}, without using quadratic mean functions. We also restrict the discretization step sizes for the hyperparameter functions $(\ell_q^{-1},\ell_p^{-1})$ to $\Delta x \in \{0,05, 0.1, 0.2, 0.25, 0.5\}$. The resulting histograms are given in Figure \ref{fig:histo_no_mean}.

\begin{figure}[h!]
    \centering
    \includegraphics[width=0.4\linewidth]{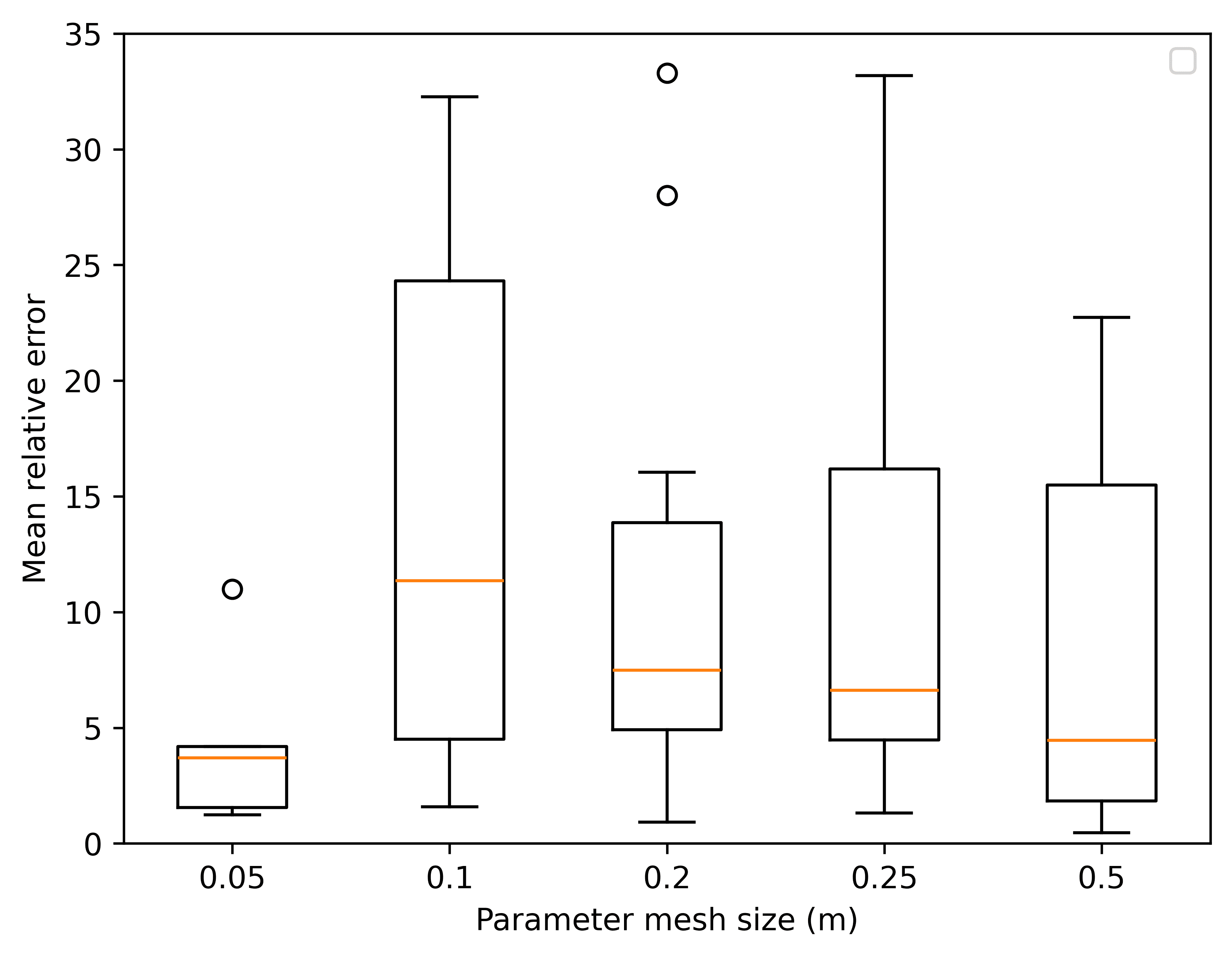}
    \caption{Study of the behaviour of PFEM-GP-dPHS as a function of the hyperparameter mesh size, without using a quadratic mean.}
    \label{fig:histo_no_mean}
\end{figure}

Numerical experiments show that without using a quadratic mean, the NLML minimisation step rarely converges to a satisfying solution $(\ell_q^{-1},\ell_p^{-1}, \sigma_f, \sigma_{\text{noise}})$. Quantitatively, we can see in Figure \ref{fig:histo_no_mean} that, even for $\Delta x = 0.05$, we often obtain a time-averaged relative $L^2$ error larger than $1$. This provides evidence that, given a fixed learning set, using a quadratic mean greatly increases the approximation quality of the Hamiltonian and more importantly, of the Hamiltonian dynamics (equations \eqref{eq:dphs-balance}). In fact, given the high dimensionality of the discretized dPHS from Section \ref{sec:num} ($d = N_p + N_q = 42$ or $82$) and the relatively high number of training points (840), we expect that from a computational perspective, obtaining a reliable discretized dPHS may become close to intractable without using a quadratic mean. To finish, this issue will become even more prohibitive for dPHS with multiple space dimensions.

\paragraph{Enforcing the boundedness of the solution $\balpha$ in the phase space}
On a theoretical level, ensuring that $\mathcal{H}(\balpha)\to + \infty$ (and that $\mathcal{H}$ is bounded from below) when $\|\balpha\| \to +\infty$ ensures that the sublevel sets $L_s \coloneqq \{ \balpha : \mathcal{H}(\balpha) \leq s\}$ are bounded for all $s\in\R$. This implies that the solutions to the associated Hamiltonian dynamics $\dot{\balpha} = (\mathcal{J}-\mathcal{R})(\delta_{\balpha} \mathcal{H}) (\balpha)$ remain bounded in phase space, \textcolor{black}{if the dPHS is well-posed in suitable function spaces. More precisely,
integrating equation \eqref{eq:dphs-balance}, using that $\mathcal{R}$ is nonnegative definite and using the Cauchy-Schwarz inequality,}
\begin{align*}
    \mathcal{H}(\balpha(t)) &= \mathcal{H}(\balpha_0) + \int_0^t \bigg(- \int_{\cal O} \bm{e}(\balpha(s)).\mathcal{R}\bm{e}(\balpha(s)) + \langle u_\partial(s), y_\partial(s) \rangle_{\partial \cal O} \bigg)\mathrm{d} s \\
    &\leq \mathcal{H}(\balpha_0) + \int_0^t \langle u_\partial(s), y_\partial(s) \rangle_{\partial \cal O} \mathrm{d}s \leq \underbrace{\mathcal{H}(\balpha_0) + \int_0^t \| u_\partial(s)\|_{L^2(\partial \mathcal{O)}} \|y_\partial(s)\|_{L^2(\partial \mathcal{O)}} \mathrm{d}s}_{=: C(\balpha_0,u_\partial,t)}.
\end{align*}
\textcolor{black}{Above, $\balpha_0$ is the initial data and $u_\partial$ is the control. If the dPHS is well-posed in suitable function spaces then $\|y(s)\|_{L^2(\partial\mathcal{O})}$ is a continuous function of $s, \balpha_0$ and $u_\partial$, and $(\balpha_0, u_\partial,t) \mapsto C(\balpha_0,u_\partial,t)$ is a continuous function.
Thus, $\balpha(t) \in L_{C(\balpha_0,u_\partial,t)}$ which is a bounded set in the phase space whose growth in size is determined by the growth of $C(\balpha_0,u_\partial,t)$.
Ensuring this boundedness property at the level of the GP prior over $\mathcal{H}$ can be understood as a form of physics-informed GP model (although not all (d)PHS enjoy this boundedness property, such as the Schrödinger equation with a finite potential well)}. In any case, this boundedness property is usually not ensured when using standard kernels alone. For example, if the kernel $k$ is of the form $k(\balpha,\balpha') = k_S(\|\balpha-\balpha'\|)$ where $k_S(h) \to 0$ when $h\to+\infty$, then regressors built from it will typically be of the form $\balpha \mapsto \sum_{k=1}^N a_jk_S(\|\balpha-\balpha_j\|)$, which also goes to $0$ when $\|\balpha\| \to + \infty$. \textcolor{black}{This implies that the set $\{ \balpha : \mathcal{H}(\balpha) \leq s\}$ is unbounded for all $s>0$, potentially allowing the solution $\balpha$ to go to infinity for new sets of inputs $(\balpha_0,u_\partial)$.}

To ensure the boundedness property above without using a quadratic mean, an alternative would be to use a non-stationary kernel e.g. of the form $k(\balpha,\balpha') = \sigma(\|\balpha\|)\sigma(\|\balpha'\|)k_S(\|\balpha-\balpha'\|)$ where $\sigma$ is a preselected function such that, for instance, $\sigma(h)k_S(h) \sim_{+\infty}  Ch^\beta$ for some chosen $C,\beta>0$. Still, this model is less ``physics-informed'' than the one with a quadratic mean, following the description given in Section \ref{subsec:select_prior_h}.
%has the disadvantage of losing the physics-informed nature of the quadratic mean.
%Even if $k$ is not stationary, the approximant above will at best remain bounded since, from the Cauchy-Schwarz inequality in the RKHS associated to $k$,
%\begin{align*}
%    \bigg|\sum_{k=1}^N a_jk(\balpha,\balpha_k)\bigg| \leq \sum_{k=1}^N |a_j|k(\balpha_j,\balpha_j)
%\end{align*}
%and thus.

%Théorie : stabilisation des trajectoires (assurer leur caractère borné) à l'aide d'un puit de potentiel quadratique : par construction, la trajectoire ne peut pas partir à l'infini.

\section{PFEM-GP-dPHS results using polynomial approximations for the GP functional hyperparameters}\label{app:polyn}
We consider the same test case as the one described in Section \ref{sec:num}. We now seek the functions $m_q,\ m_p,\ \ell_q^{-1}$ and $\ell_p^{-1}$ in the space of polynomial functions of order 3 instead of the space of piecewise linear functions. This results in $4\times 4 + 2 =18$ hyperparameters in total, as opposed to e.g. 26 hyperparameters estimated to obtain Figure~\ref{fig:hyperparameters}. The obtained estimations are given in Figure~\ref{fig:gp_polyn}, which is a counterpart of Figure~\ref{fig:hyperparameters}.
\begin{center}
    \begin{figure}[h!]
    \includegraphics[width=0.45\textwidth]{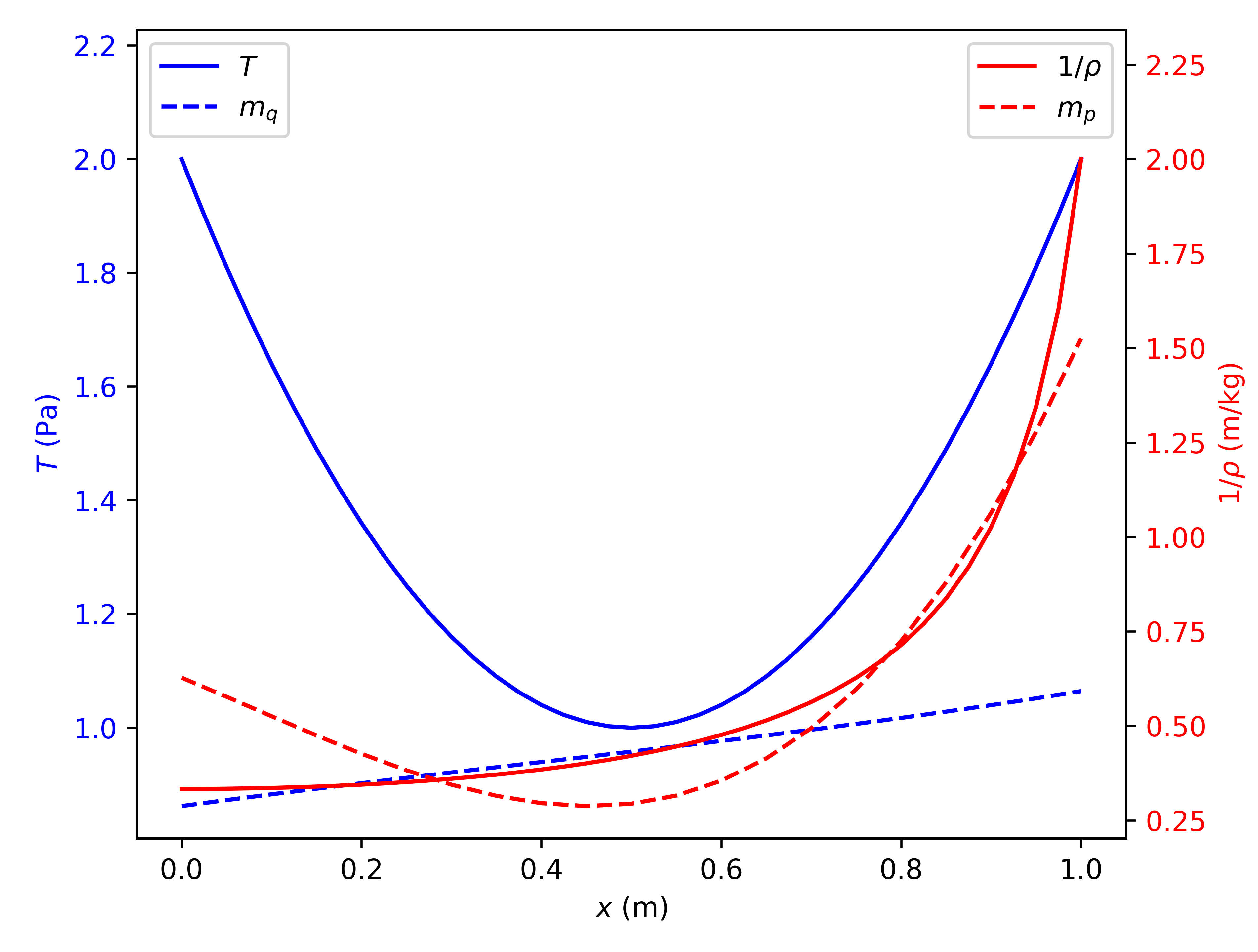}\centering
    \includegraphics[width=0.45\textwidth]{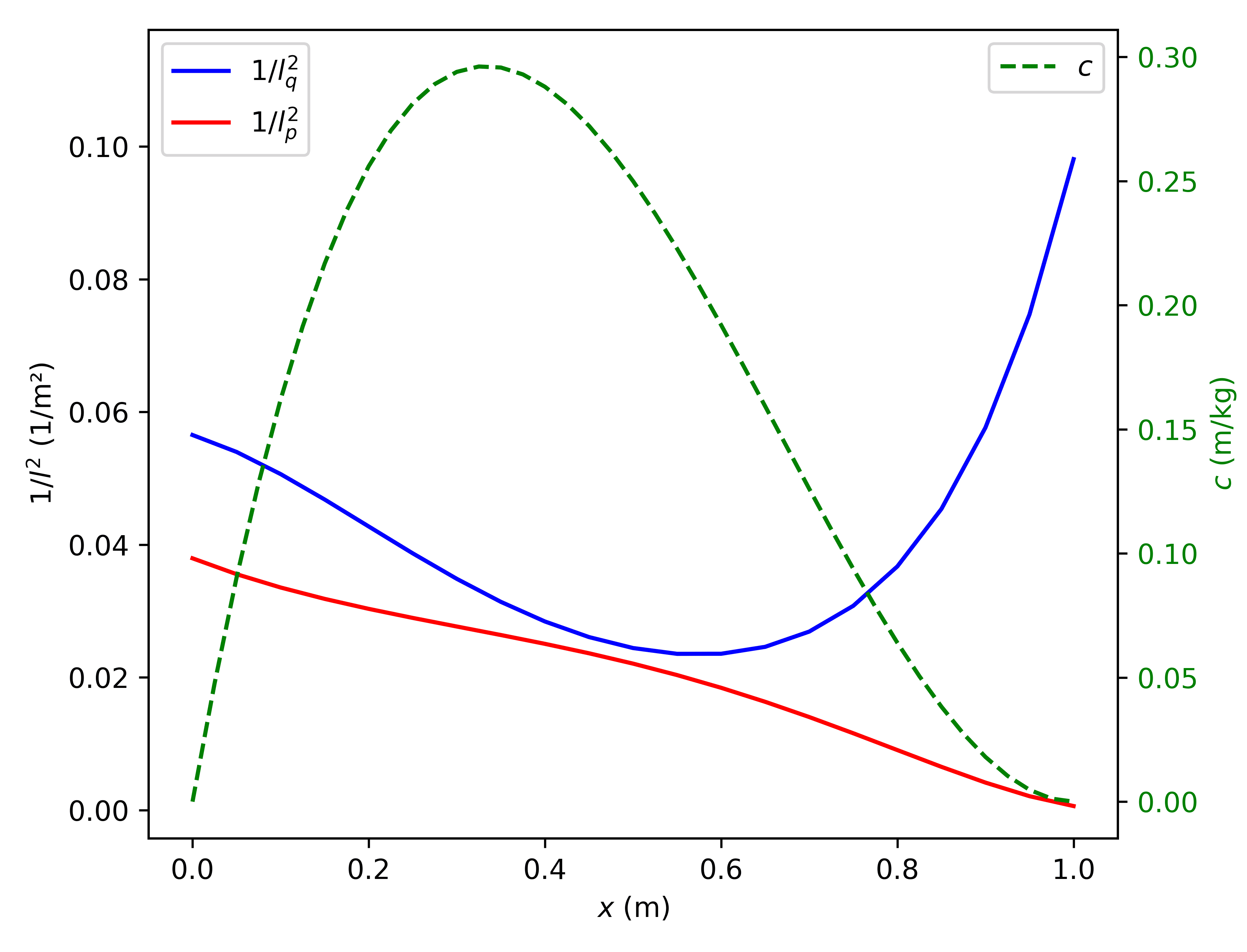}
    \caption{Left : estimation of the GP quadratic mean functions $m_q$ and $m_p$ using polynomial functions of order 3, versus the true physical parameters $T$ and $\rho^{-1}$. Right : estimation of the GP lengthscale functions $\ell_q^{-1}$ and $\ell_p^{-1}$ using  polynomial functions of order 3, versus the function $c$.}
    \label{fig:gp_polyn}
\end{figure}
\end{center}
Figure~\ref{fig:gp_polyn} (left) shows that $\rho^{-1}$ is more or less well-estimated by $m_p$. However, the function $T$ is not well approximated by $m_q$. As a result, the quadratic part of the Hamiltonian is not identified correctly by the prior mean with tuned hyperparameters. Contrarily to the piecewise linear case, this implies that the kernel $k$ must be active even at locations where the quadratic part of the Hamiltonian is dominant. In particular, there is no reason for the correlation lengths to follow the variations of $c$, as shown in Figure \ref{fig:gp_polyn} (right). Yet, Figure~\ref{fig:alpha_p_mix} (right) provides visual evidence that, even in this case where the kernel needs to compensate for an ill-estimated quadratic mean, the PFEM-GP-dPHS framework is still capable of providing a good estimation of ~$\alpha_p$.

Possible extensions of this study include using other bases, such as the Fourier basis or wavelets, and studying their performance as a function of their discretization mesh size.

\section{Extension to 2D dPHS}\label{app:2d}

Extension of the PFEM-GP-dPHS method to the 2D case, and generally higher dimension case, is possible using the \texttt{SCRIMP} package. We show here the calculation and numerical result for the nonlinear 2D wave equation, following the same workflow as the 1D case in Sections \ref{sec:discretization-1D}, \ref{sec:prior-pfem-gp-dphs} and \ref{sec:num}.

\subsection{Nonlinear 2D wave equation in dPHS form}\label{subsec:2Ddphs}

We consider the direct extension our 1D exemple, the 2D wave equation with Neumann boundary controls on a bounded rectangle $\O_{2D} := (0, L) \times (0, \ell)$. The deflection of the membrane from the equilibrium $w$ satisfies 
\begin{equation}
    \left\lbrace
    \begin{array}{rcl}
    \rho(x)\,\partial_{tt}^2 w(t,x) - \mathrm{div} \left( \overline{\overline{h}}(x, {\mathbf{grad}}  \left( w(t,x) \right) ) \cdot {\mathbf{grad}}  \left( w(t,x) \right) \right) = 0, \ \ \ t \ge 0, x \in \O_{2D}, \\
    \bm{n} \cdot \overline{\overline{h}}(s, {\mathbf{grad}}  \left( w(t,s) \right) ) \cdot {\mathbf{grad}}  \left( w(t,s) \right) = u_N(t,s), \ \ \  t \ge 0, s \in \Gamma = \delta\O_{2D} , \\
    \end{array}
    \right.
\end{equation}
We set the energy variable $\balpha_q := {\mathbf{grad}}  w \in \R^2$ the strain, and $\alpha_p := \rho \partial_t w \in \R$ the linear momentum. Note that the main difference with the 1D case is that $\balpha_q$ is now a vector-valued variable.

The term $\overline{\overline{h}}(x, \balpha_q(t,x)  ) = \overline{\overline{T}}(x) + \overline{\overline{c}}(x)\,g(\balpha_q(t,x))$  is made of a linear term (Elasticity tensor) $\overline{\overline{T}}(x) \in \R^{2\times 2}$, and a nonlinear term $\overline{\overline{c}}(x)\,g(\balpha_q(.,x))$ with $\overline{\overline{c}}(x) \in \R^{2\times2}$ and $g(\bm{u}) = \exp(-\bm{u}\cdot \bm{u})$.  One can express the total mechanical energy $\mathcal{H}$, the Hamiltonian, as
\begin{equation}
    \mathcal{H}(t) = \mathcal{H}(\balpha_q(t,x), \alpha_p(t,x)) := \underbrace{\frac{1}{2} \int_{\O_{2D}} \frac{\alpha_p(t,x)^2}{\rho(x)} \mathrm{d}x}_{\text{Kinetic energy}} + \underbrace{\frac{1}{2} \int_{\O_{2D}} \balpha_q^\top(t,x)  \overline{\overline{T}}(x)  \balpha_q(t,x) \mathrm{d}x}_{\text{Quadratic potential energy}} \nonumber
\end{equation}

\begin{equation}
+\underbrace{ \frac{1}{2} \int_{\O_{2D}} \int_{0<y<\alpha_{q1}(t,x), \,0<z<\alpha_{q2}(t,x)} g(\begin{bmatrix}
    y \\ z
\end{bmatrix})\begin{bmatrix}
    y & z
\end{bmatrix}\overline{\overline{c}}(x) \begin{bmatrix}
    \mathrm{d}y \\ \mathrm{d}z
\end{bmatrix}}_{\text{Nonlinear term}} .
\end{equation}
The equation above is a 2D generalization of equation \eqref{eq:hamiltonian_simu}.
The co-energy variables are, as in the 1D case $\bm{e}_q := \delta_{\balpha_q} \mathcal{H} = \overline{\overline{h}} \cdot \balpha_q $ the stress, and $ e_p := \delta_{\alpha_p} \mathcal{H} = \frac{\alpha_p}{\rho}$ the velocity.

Using Newton's second law and Schwarz's lemma, we obtain port-Hamiltonian system representing a vibrating membrane with Neumann boundary controls :

\begin{equation}
    \begin{pmatrix} \partial_t \balpha_q \\ \partial_t \alpha_p \end{pmatrix}
    =
    \begin{bmatrix} 0 & \mathbf{grad} \\ \mathrm{div} & 0 \end{bmatrix}
    \begin{pmatrix}\bm{e}_q \\ e_p \end{pmatrix},
\qquad
    \left\lbrace
    \begin{array}{rcl}
   \bm{e}_q(t,s)\cdot\bm{n} &=& u(t,s), \qquad t \ge 0, \,s \in \Gamma, \\
    y(t,s) &=& e_p(t,s), \qquad t \ge 0, \,s \in \Gamma. \\
    \end{array}
    \right.
    \label{eq:dPHS_2D}
\end{equation}

The power balance satisfied by the Hamiltonian is
\begin{equation}
    \frac{\mathrm{d}}{\mathrm{d}t} \mathcal{H}(t) = \frac{\mathrm{d}}{\mathrm{d}t} \mathcal{H}(\balpha_q(t), \alpha_p(t)) = \underbrace{\left\langle y(t,\cdot), u(t,\cdot)\right\rangle_{\Gamma}}_{\text{power flowing through }\Gamma}.
\end{equation}

\subsection{PFEM}\label{subsec:2DPFEM}

Now let apply the PFEM on this 2D problem. Let $\bvarphi_q$ and $\varphi_p$ be smooth test functions on $\O_{2D}$, and $\psi$ be smooth test functions on $\Gamma$ . One can write the weak formulation of the Dirac structure as follows :
\begin{equation}
    \left\lbrace
    \begin{array}{rcl}
    \int_{\O_{2D}} \partial_t \balpha_q(t,x) \cdot \bvarphi_q(x) \mathrm{d}x &=& \int_{\O_{2D}} {\mathbf{grad}}  \left( e_p(t,x) \right) \cdot \bvarphi_q(x) \mathrm{d}x, \\
    \int_{\O_{2D}} \partial_t \alpha_p(t,x) \varphi_p(x) \mathrm{d}x &=& \int_{\O_{2D}} \mathrm{div} \left(\bm{e}_q(t,x) \right) \varphi_p(x) \mathrm{d}x, \\
    \left\langle y, \psi \right\rangle_{\Gamma} &=& \left\langle e_p, \psi \right\rangle_{\Gamma}\,.
    \end{array}
    \right.
\end{equation}

Integrating the second line  by parts makes the control $u$ and the observation $y$ appear:
\begin{equation}
    \int_{\O_{2D}} \partial_t \alpha_p(t,x) \varphi_p(x) \mathrm{d}x = - \int_{\O_{2D}}\bm{e}_q(t,x) \cdot {\mathbf{grad}}  \left( \varphi_p(x) \right) \mathrm{d}x + \left\langle \varphi_p, u \right\rangle_{\Gamma} 
\end{equation}

Let $(\bvarphi_q^i)_{1 \le i \le N_q} \subset \left(L^2(\O_{2D})\right)^2$ and $(\varphi_p^k)_{1 \le k \le N_p} \subset H^1(\O_{2D})$ be two finite families of approximations so that $ \balpha_q(t,x) \simeq  \underline{\varphi_q }^\top(x) \underline{\alpha_q}(t)$, $\alpha_p(t,x) \simeq \underline{\varphi_p }^\top(x) \underline{\alpha_p}(t)$ (likewise for $\dot{\balpha_q}$ and $\dot{\alpha_p}$) and $\bm{e}_q(t,x) \simeq  \underline{\varphi_q }^\top(x) \underline{e_q}(t)$, $e_p(t,x) \simeq \underline{\varphi_p }^\top(x) \underline{e_p}(t)$. One thing to remark is that compared to the 1D case, $\underline{\varphi_q }^\top = [\bvarphi_q^1, \cdots,\bvarphi_q^{N_q}] \in \left(L^2(\O_{2D})\right)^{2 \times N_q} $ is now a matrix. We denote also $(\psi^m)_{1 \le m \le N} \subset H^{\frac12}(\Gamma)$. In particular, the latter choices imply that the duality brackets at the boundary reduce to simple $L^2$ scalar products. We next replace the variables by their approximations, and write the system in matrix form

\begin{equation}
    \left\lbrace
    \begin{array}{rcl}
    \underbrace{\begin{bmatrix}
    M_q & 0  \\
    0 & M_p 
    \end{bmatrix}}_{M} 
     \underbrace{\begin{pmatrix}
    \underline{\dot{\alpha}_q}(t) \\
    \underline{\dot{\alpha}_p}(t)
    \end{pmatrix}}_{ \underline{\dot{\alpha}}}
    &=&
   \underbrace{ \begin{bmatrix}
    0 & D &  \\
    -D^\top & 0
    \end{bmatrix}}_{J-R}
     \underbrace{ \begin{pmatrix}
    \underline{e_q}(t) \\
    \underline{e_p}(t) 
    \end{pmatrix}}_{ \underline{e}} + \underbrace{ \begin{bmatrix}
    0  \\
    B 
    \end{bmatrix}}_{G}
    \underline{u}, \\
   \underline{y}(t) 
    &=&
   \begin{bmatrix}
   
    0 & B^\top \\
   
    \end{bmatrix}
    \begin{pmatrix}
    \underline{e_q}(t) \\
    \underline{e_p}(t) 
    \end{pmatrix},
    \end{array}
    \right.
\end{equation}
with square matrices $ (M_q)_{ij} := \int_{\O_{2D}} \bvarphi_q^j(x)\cdot \bvarphi_q^i(x) {\rmd}x$,~ $ (M_p)_{k\ell} := \int_{\O_{2D}} \varphi_p^\ell(x) \varphi_p^k(x) {\rmd}x$ and $ (M_\Gamma)_{k\ell} := \int_\Gamma \psi^\ell(s) \psi^k(s) {\rmd}s$, and rectangular matrices~ $(D)_{i\ell} := \int_{\O_{2D}} {\mathbf{grad}} \left( \varphi_p^\ell(x) \right) \cdot\bvarphi_q^i(x) {\rmd}x$ and $(B)_{k\ell} := \int_\Gamma \varphi_p^\ell(s) \psi^k(s) {\rmd}s$.
The system is closed by approximating the constitutive relations like in 1D by writing them in their weak form:

\begin{equation}
    \left\lbrace
    \begin{array}{rcl}
    \int_{\O_{2D}}\bm{e}_q(t,x) \cdot \bvarphi_q(x) {\rmd}x &=& \int_{\O_{2D}} \left(\overline{\overline{h}}(x,\balpha_q(t,x)) \,\balpha_q(t,x)\right) \cdot \bvarphi_q(x) {\rmd}x, \\
    \int_{\O_{2D}} e_p(t,x) \varphi_p(x) {\rmd}x &=&  \int_{\O_{2D}} \frac{1}{\rho(x)} \alpha_p(t,x) \varphi_p(x) {\rmd}x.
    \end{array}
    \right.
\end{equation}

\noindent
Then the matrix form of the discrete weak formulation of the constitutive relations is

\begin{equation}
       \left\lbrace
    \begin{array}{rcl}
    M_q \underline{e_q}(t) &=& M_h \underline{\alpha_q}(t), \\
    M_p \underline{e_p}(t) &=& M_\rho \underline{\alpha_p}(t),
    \end{array}
    \right.
\end{equation}
 
\noindent
where $(M_h)_{ij} := \int_{\O_{2D}}  \bvarphi_q^j(x)\cdot \overline{\overline{h}}(x,\balpha_q(t,x)) \cdot \bvarphi_q^i(x) {\rmd}x$ ~and~ $(M_\rho)_{k\ell} := \int_{\O_{2D}} \frac{1}{\rho(x)} \varphi_p^\ell(x) \varphi_p^k(x) {\rmd}x$.

\subsection{PFEM-GP-dPHS}\label{subsec:PFEM-GP_2D}

To write the associated PFEM-GP-dPHS as in Section ~\ref{sec:prior-pfem-gp-dphs}, we use a scalar prior mean $m$ and kernel $k$ selected for the GP of $\widehat{\mathcal{H}}$, such that $\widehat{\mathcal{H}} \sim \GP(m,k)$, with

\begin{align} 
m(\balpha) &=  \int_{\O_{2D}}\balpha_q^\top(x) A_{q}(x)\balpha_q(x)\,dx +  \int_{\O_{2D}} m_{p}(x)\alpha_p^2(x)\,dx, 
\label{eq:mean_GP_2D}\\
    k(\balpha,\balpha') &= \sigma_f^2 \exp\bigg[-\frac{1}{2}\bigg( \int_{\O_{2D}} (\balpha_q(x)-{\balpha_q'}(x))^\top L_q^{-2}(x) (\balpha_q(x)-{\balpha_q'}(x))\,dx \nonumber \\
    & \hspace{6.5cm} +  \int_{\O_{2D}} \frac{(\alpha_p(x)-{\alpha_p'}(x))^2)}{\ell_p^2(x)}\,dx \bigg)\bigg], 
\end{align}
where the functions $m_{p}$ and $\ell_p$ are 2D nonnegative hyperparameters, and the PSD-matrix-valued function $A_{q}(x) = \left[[m_{q1}(x) ~m_{q3}(x)]^\top~~[m_{q3}(x) ~ m_{q2}(x)]^\top\right]
$ and diagonal matrix-valued function $L_q(x) = \text{Diag}(\ell_{q1}(x), \ell_{q2}(x))\in\R^{2\times 2}$ are matrix hyperparameters of the GP $\widehat{\mathcal{H}}$.

Next, we adapt the affine transformation of equation ~\eqref{eq:transfo} we need to obtain the GP distribution of $e$:

\begin{equation}
    M\underline{e}(\balpha) = \int_{\O_{2D}}\underbrace{\begin{bmatrix} \underline{\varphi_q}(x) &  0_{N_q,1}\\ 0_{N_p,2} & \underline{\varphi_p}(x) \end{bmatrix} }_{\Phi(x)} (\delta_{\balpha} \widehat{\mathcal{H}})(\balpha)(x)  {\rmd}x ,
    \label{eq:transfo_2D}
\end{equation}
where now $\Phi \in \R^{(N_q+N_p)\times3}$. We apply the variational derivative $\delta_{\balpha} = [\delta_{\balpha_q},\delta_{\alpha_p}]^\top$ to obtain the prior GP distribution over $e$

\begin{equation}
    m_e(\balpha) = \begin{bmatrix}A_{q}\balpha_q \\ m_{p}\alpha_p  \end{bmatrix}\,,
\end{equation}

\begin{align}
    \frac{k_e(\balpha,\balpha')}{ k(\balpha,\balpha')}  &=   \begin{bmatrix}
 \text{Diag}\Big( \frac{1}{\ell_{q1}^2}\mathcal{I}, \frac{1}{\ell_{q2}^2}\mathcal{I}\Big)   & 0 \\
   0& \frac{1}{\ell_p^2}\mathcal{I} 
\end{bmatrix}\\
& \ \ \ \ \ -\begin{bmatrix}
 L_q^{-2}(\balpha_q -\balpha_q') \otimes L_q^{-2}(\balpha_q-\balpha_q')^\top  & L_q^{-2}(\balpha_q -\balpha_q') \otimes \frac{1}{\ell_p^2}(\alpha_p-\alpha_p') \\
    \frac{1}{\ell_p^2}(\alpha_p -\alpha_p') \otimes L_q^{-2}(\balpha_q-\balpha_q')^\top&  \frac{1}{\ell_p^2}(\alpha_p -\alpha_p') \otimes \frac{1}{\ell_p^2}(\alpha_p-\alpha_p') 
\end{bmatrix}
\end{align}

Finally, we multiply by $\Phi$, integrate and replace $\balpha$ by its finite-dimensional approximation $\Phi^\top\underline{\alpha}$ to obtain the prior over $\underline{e}(\underline{\alpha})$

\begin{equation}
   m_{\underline{e}}(\underline{\alpha}) = \begin{bmatrix} (M_{m,q} \underline{\alpha_q})^\top &  (M_{m,p} \underline{\alpha_p})^\top   \end{bmatrix}^\top,
   \label{eq:mean_2d}
\end{equation}

\noindent with $(M_{m,q})_{ij} :=  \int_{\O_{2D}}\bvarphi_q^{i\top}  A_{q}\bvarphi_q^j {\rmd}x$ and $(M_{m,p})_{ij} :=  \int_{\O_{2D}} m_{p} \varphi_p^i \varphi_p^j{\rmd}x$ and 

\begin{equation}
        \frac{k_{\underline{e}}(\underline{\alpha}, \underline{\alpha}')}{k(\Phi^\top\underline{\alpha},\Phi^\top\underline{\alpha}')} =  M_k\begin{bmatrix}
       I_{N_q}-(\underline{\alpha_q}-\underline{\alpha_q'})(\underline{\alpha_q}-\underline{\alpha_q'})^\top  &-(\underline{\alpha_q}-\underline{\alpha_q'})(\underline{\alpha_p}-\underline{\alpha_p'})  \\-(\underline{\alpha_p}-\underline{\alpha_p'})(\underline{\alpha_q}-\underline{\alpha_q'})^\top &I_{N_p}-(\underline{\alpha_p}-\underline{\alpha_p'})(\underline{\alpha_p}-\underline{\alpha_p'})^\top
       \end{bmatrix}M_k^\top,
       \label{eq:kernel_2d}
\end{equation}

\noindent where  $M_k = M^{-1}\text{Diag}(M_{k,q},M_{k,p})$, where the submatrices $M_{k,q}$ and $M_{k,p}$ are given by $(M_{k,q})_{ij} :=  \int_{\O_{2D}}  \bvarphi_q^{i\top}L_q^{-2} \bvarphi_q^j{\rmd}x$ and $(M_{k,p})_{ij} :=  \int_{\O_{2D}} \ell_p^{-2} \varphi_p^i \varphi_p^j{\rmd}x$.

Finally to obtain the GP prior over $ f_{\dot{\underline{\alpha}}}$, like in 1D we  apply the matrix multiplication $(J-R)$ to the GP $\underline{e}$ and add the mean $G\underline{u}$. 

\subsection{Generalization to higher dimensions}\label{subsec:PFEM_GP_nD}
From here, we can remark that the only difference between the 1D PFEM-GP-dPHS and 2D PFEM-GP-dPHS lies only in the matrix $M_{m,q}$,$M_{m,p}$,$M_{k,q}$ and $M_{k,p}$, as equations ~\eqref{eq:mean_e} and ~\eqref{eq:mean_2d}, and also equations ~\eqref{eq:kernel_e} and ~\eqref{eq:kernel_2d}, are the same. From here, it is easy to see than for any dPHS over a space $\O_{nD}$ of dimension $n$D and with vector-valued energy variables of size $n_q$ and $n_p$, the associated PFEM-GP-dPHS will use the mean and kernel functions given in equations ~\eqref{eq:mean_e} and ~\eqref{eq:kernel_e}, but with 

$$(M_{m,q})_{ij} :=  \int_{\O_{nD}}\bvarphi_q^{i\top}(x)  A_{q}(x)\bvarphi_q^j(x) {\rmd}x,$$
$$(M_{m,p})_{ij} :=  \int_{\O_{nD}} \bvarphi_p^{i\top}(x)  A_{p}(x)\bvarphi_p^j(x){\rmd}x,$$
$$(M_{k,q})_{ij} :=  \int_{\O_{nD}}  \bvarphi_q^{i\top}(x)L_q^{-2} (x)\bvarphi_q^j(x){\rmd}x,$$
$$(M_{k,p})_{ij} :=  \int_{\O_{nD}} \bvarphi_p^{i\top}(x)L_p^{-2}(x) \bvarphi_p^j(x){\rmd}x,$$
where $ A_{q}$ and $ A_{p}$ are PSD-symmetric matrix-valued function of size respectively $n_q$ and $n_p$ and   $ L_{q}$ and $ L_{p}$ are diagonal matrix-valued function of size respectively $n_q$ and $n_p$.  From here, those new hyperparameter functions can, like in 1D, be discretized using their own finite element method, or even other kind of decomposition.

\subsection{Numerical results}\label{subsec:2Dnum}

We test PFEM-GP-dPHS on the dPHS described by equation \eqref{eq:dPHS_2D}, using the same numerical workflow as for the 1D equation in Section \ref{sec:num}. The space $\O_{2D}$ is meshed with 22 elements. $\balpha_q$ is discretized with $\mathbb{P}^1$ Lagrange elements giving $N_q = 132$ degrees of freedom, and $\alpha_p$ is discretized with $\mathbb{P}^2$ Lagrange elements giving $N_q = 57$ degrees of freedom. We also discretize the boundary terms $u$ and $y$ with 1D $\mathbb{P}^1$ Lagrange elements, and use a simulation time step $\Delta t = 0.01$s. The parameters of the dPHS are taken constant with $\rho^{-1} = 0.25$, $\overline{\overline{T}} = \left[[5 ~2]^\top~~[2~7]^\top\right]$ and $\overline{\overline{c}} = \left[[2 ~1]^\top~~[1~2.5]^\top\right]$. The hyperparameters are also taken to be constant over the space, giving us as hyperparameters
$A_{q} = \left[[m_{q1} ~m_{q3}]^\top~~[m_{q3} ~ m_{q2}]^\top\right]
$, $L_q = \text{Diag}(\ell_{q1}, \ell_{q2})$, $m_{p}$, $\ell_{p}^{-1}$, $\sigma_f$ and the noise, thus 9 hyperparameters. We use \texttt{SCRIMP} to generate 10s of simulation with all controls $u$ fixed to 0 and with initial condition $\balpha_q = [0.1x(x-2)y(y-1),-0.1x(x-2)y(y-1)]^\top$ and $\alpha_p = 3^{(-20((x-0.5)(x-0.5)+(y-0.5)(y-0.5)))}$. Then we select $N_s = 55$ time steps within 0 to 5 seconds to compute our training data, resulting in $N_{train} = (N_q+N_p)\times N_s = 10395$ training points. Using the same method as for the 1D case, we present in Figure \ref{fig:2D_error} and \ref{fig:2D_screensot} the result of a successful regression. Here we obtained $A_{q} = \left[[0.609 ~0.606]^\top~~[0.606 ~ 0.741]^\top\right]
$, $L_q = \left[[2.97 ~0]^\top~~[0 ~ 2.05]^\top\right]$, $m_{p} = 0.126$, $\ell_{p} = 16.10$ and $\sigma_f = 0.879$. Figure \ref{fig:2D_error} (left) shows that the relative errors of the calculated trajectory versus the true trajectory of \texttt{SCRIMP} $ \alpha^{\text{SCP}}$, $e^{\text{rel}}_{mj} = \|\alpha_m(t_j) - \alpha_m^{\text{SCP}}(t_j)\|_{L^2}/\|\alpha_m^{\text{SCP}}(t_j)\|_{L^2}$ ($m=q,p$), slowly increase but stay below 50\% within the time interval during which we measured training data. After the 5s mark, we can see that the relative error rises faster as the trajectory reaches parts of the $\balpha$ space on which the GP has not been trained. However, this divergent behaviour stays relatively slow, which can be explained by the right part of Figure \ref{fig:2D_error}. Here, we see that the relative error of the posterior mean of the GP $\underline{e}$, $\tilde{m}_{\underline{e}}$ over the true trajectory, $e^{\text{rel,e}}_{mj} = \|\tilde{m}_{\underline{e},m}(\alpha^{\text{SCP}}(t_j)) - e_m^{\text{SCP}}(\alpha^{\text{SCP}}(t_j))\|_{L^2}/\|e_m^{\text{SCP}}(\alpha^{\text{SCP}}(t_j))\|_{L^2}$ ($m=q,p$), is low before 5 seconds as it is where the training points are located, but also not very high after 5 seconds. This prevents the trajectory from diverging too fast from the \texttt{SCRIMP} simulation for ``large'' times, and also underlines the capacity of the PFEM-GP-dPHS to approximate the Hamiltonian in places of the $\balpha$ space with no training points in it. Figure \ref{fig:2D_screensot} shows a representation of the  estimated and true energy variable at $t=1$ ($\alpha_{q1}$ and $\alpha_p$), to have a qualitative appreciation of the quality of the regression.

\begin{center}
    \begin{figure}[t!]
    \includegraphics[width=0.45\textwidth]{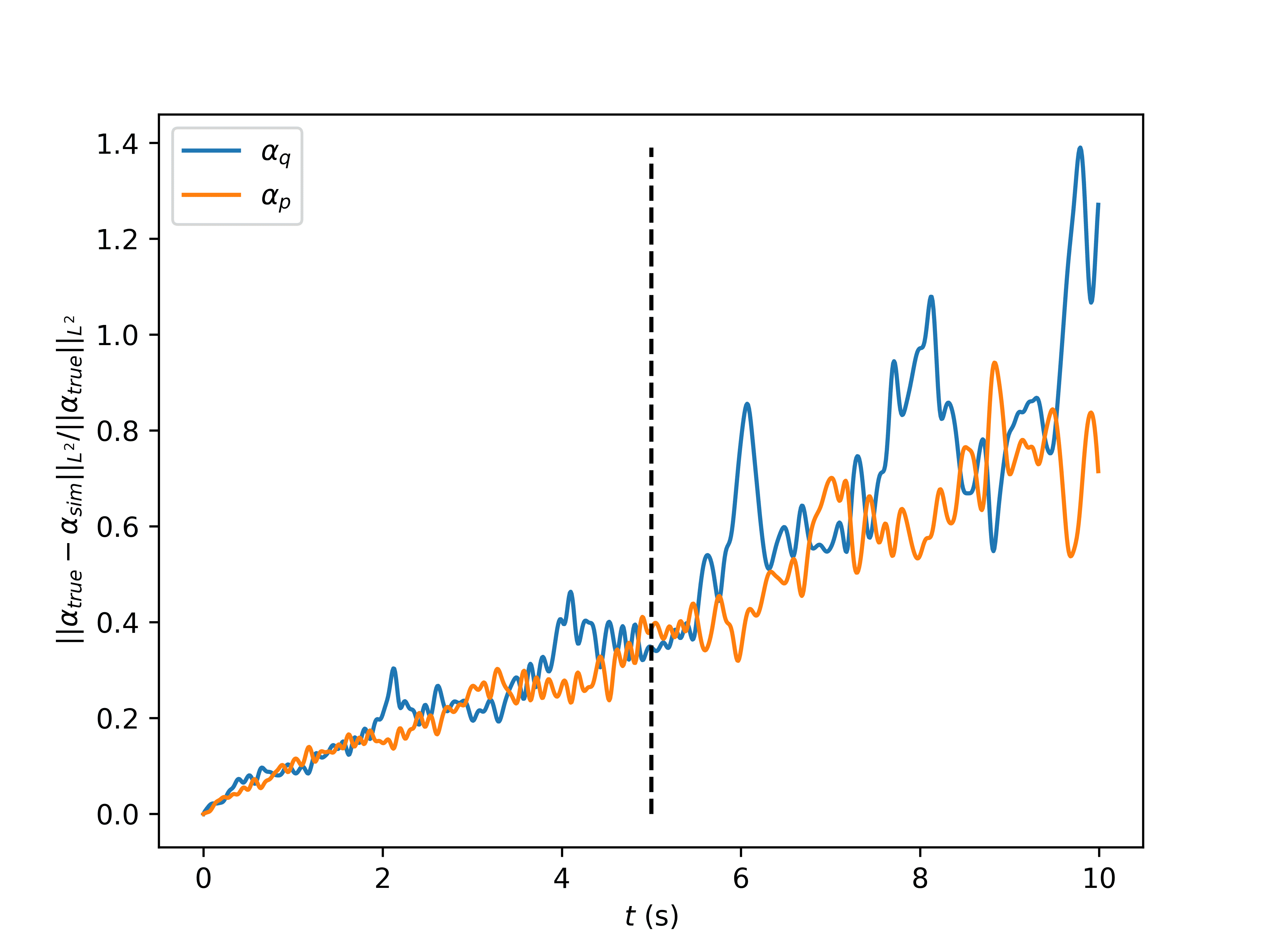}\centering
    \includegraphics[width=0.45\textwidth]{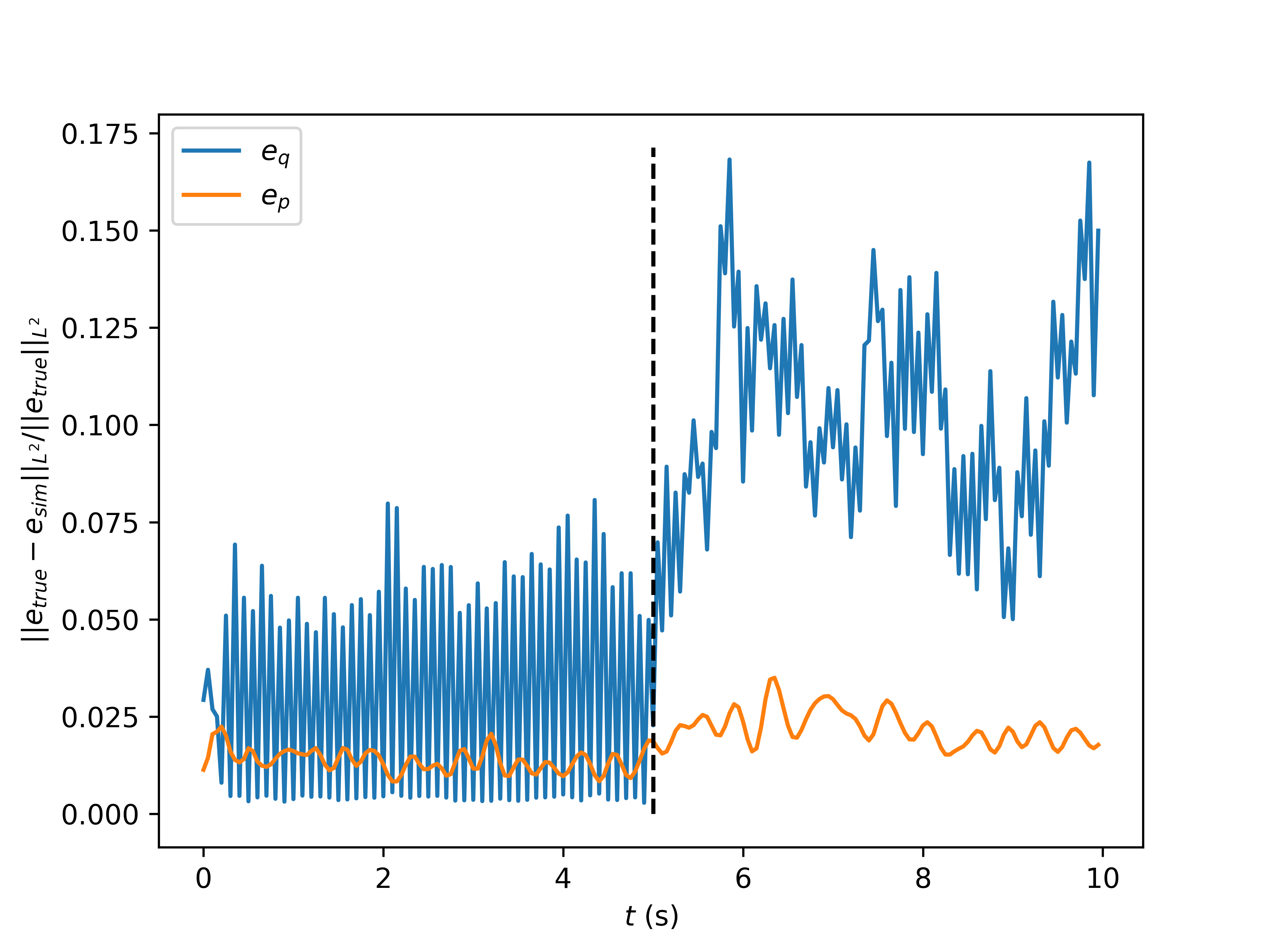}
    \caption{(Left) Evolution of the relative error of the $\balpha_q$ and $\alpha_p$ trajectories over time. (Right) Relative error of the posterior mean of $\bm{e}_q$ and $e_p$ calculated at each point of the true trajectory $\balpha (t)$. The dotted black line represents the time where we stop using time step for input data.}
    \label{fig:2D_error}
\end{figure}
\end{center}

\begin{center}
    \begin{figure}[h!]
    \includegraphics[width=0.45\textwidth]{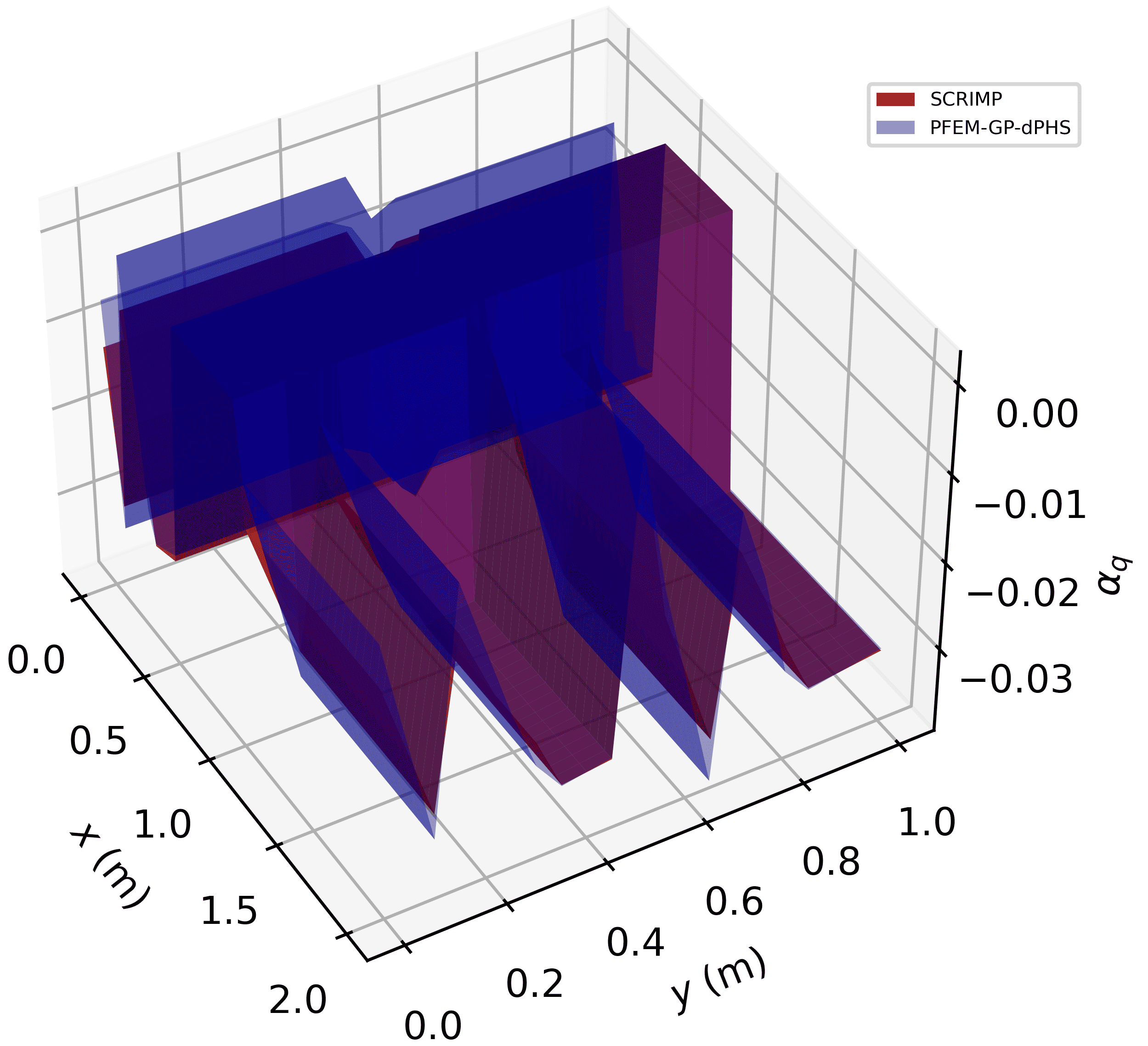}\centering
    \includegraphics[width=0.45\textwidth]{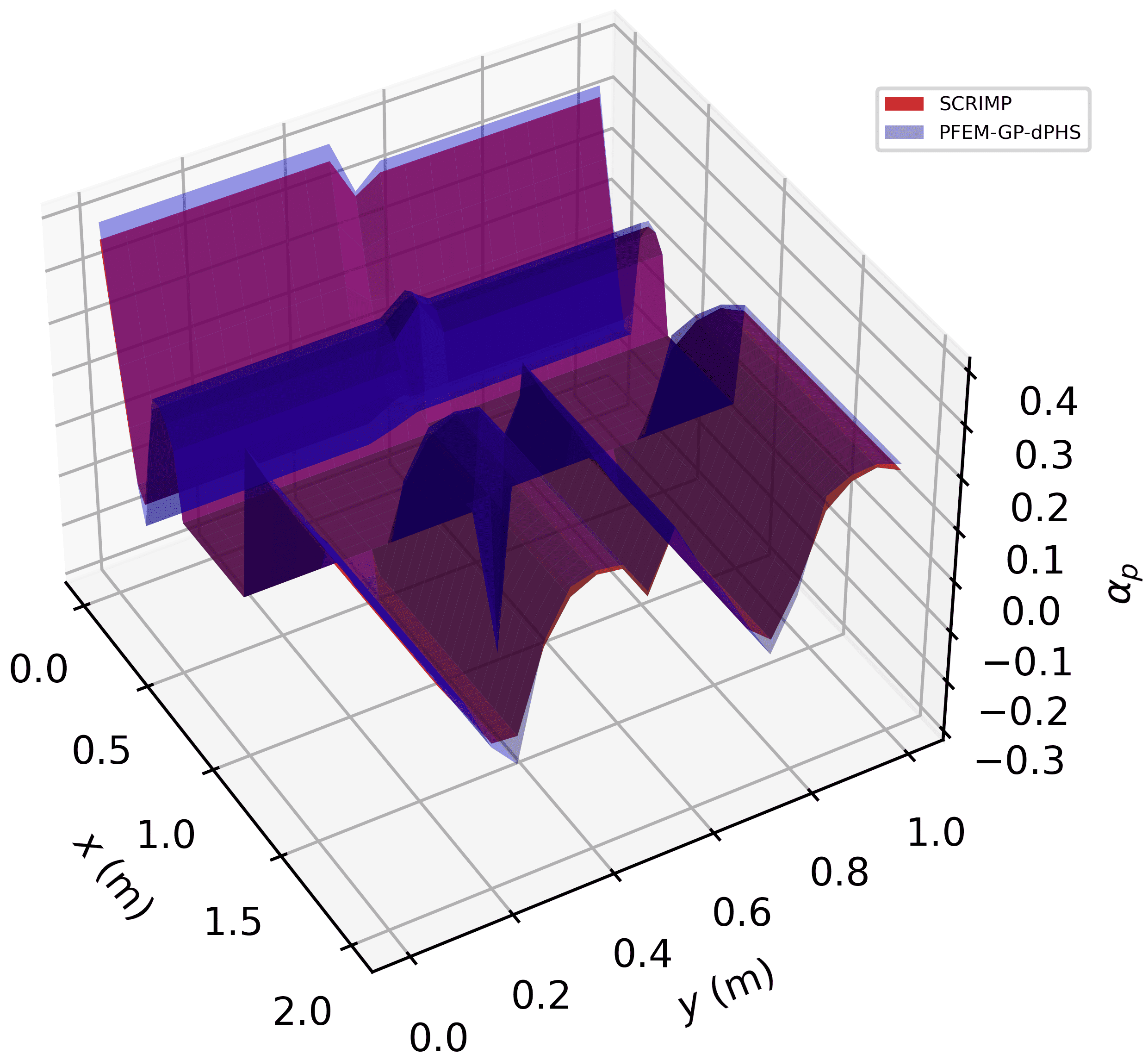}
    \caption{Representation of the true and simulated (left) $\alpha_{q1}$ and (right) $\alpha_{p}$ at time $t=1s$ for  for the 2D case.}
    \label{fig:2D_screensot}
\end{figure}
\end{center}

\section{Computational costs and additional 1D simulation figures}\label{app:additional_simu}
\subsection{Additional 1D simulation figures\label{subsec:additional_figs}}

Figures \ref{fig:alpha_p_mix} and \ref{fig:alpha_q_mix} present additional 1D simulation outputs, showing the trajectory of $\alpha_p$ or $\alpha_q$ over time for the SCRIMP simulation and the predicted trajectory after regression with the PFEM-GP-dPHS method. Figure \ref{fig:alpha_p_mix} shows two successful predicted trajectories using 2 different decomposition methods: piecewise linear functions and 3rd order polynomial functions. The analysis of this result is done in Appendix \ref{app:polyn}. Figure \ref{fig:alpha_q_mix} illustrates the behaviour of PFEM-GP-dPHS as a function of the number of finite elements: using a finer finite element mesh for the energy variables allows to obtain smoother trajectories after regression.

\begin{figure}[h!]
    \centering
    \includegraphics[width=0.9\linewidth]{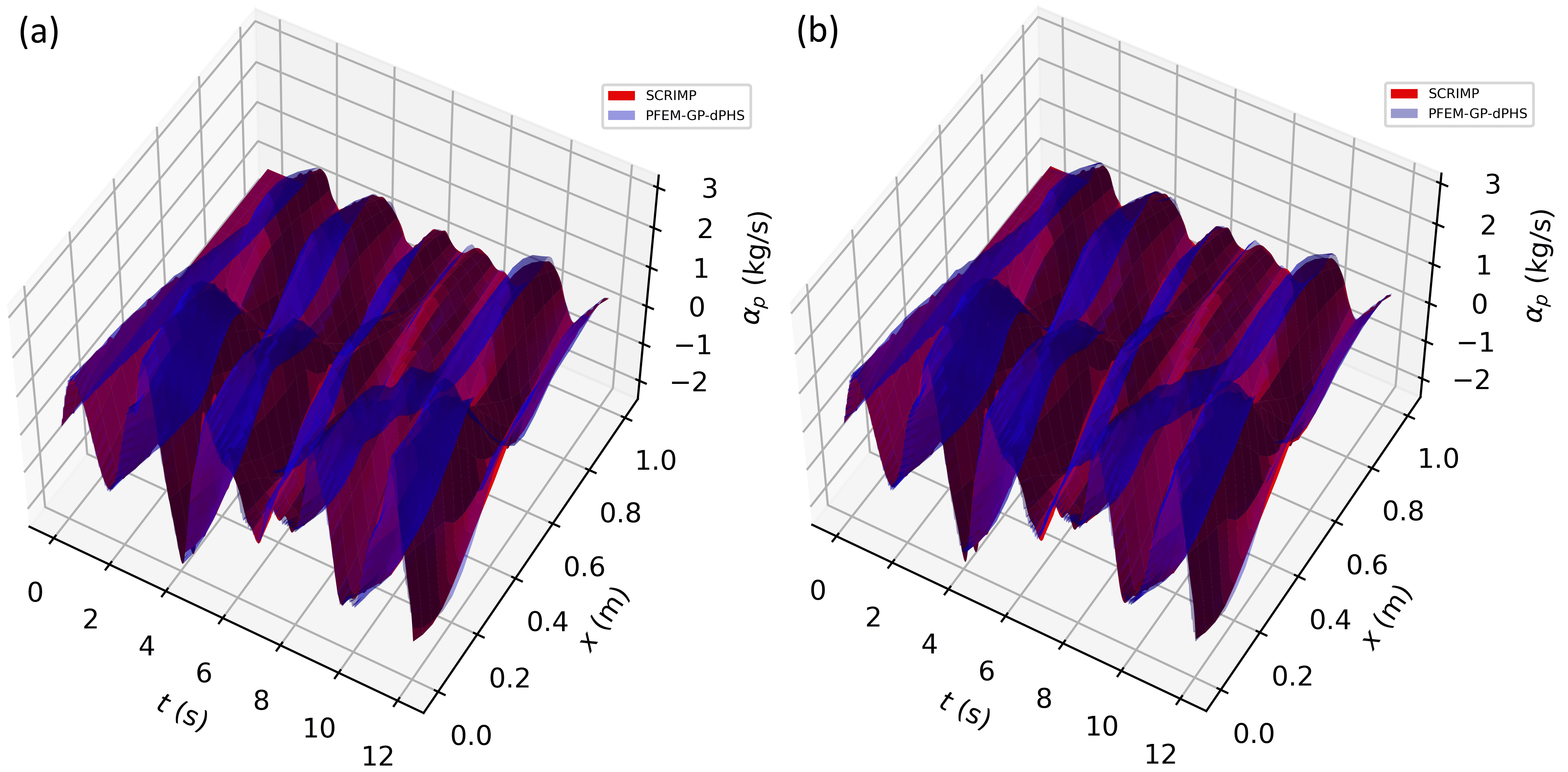}
    \caption{True $\alpha_p$ trajectory obtained with \texttt{SCRIMP} plotted with the estimated $\alpha_p$ trajectory obtained using the posterior mean of the PFEM-GP-dPHS and discretizing hyperparameters with (left) piecewise linear functions (right) 3 order polynomial functions.}
    \label{fig:alpha_p_mix}
\end{figure}

\begin{center}
    \begin{figure}[h!]
    \includegraphics[width=0.45\textwidth]{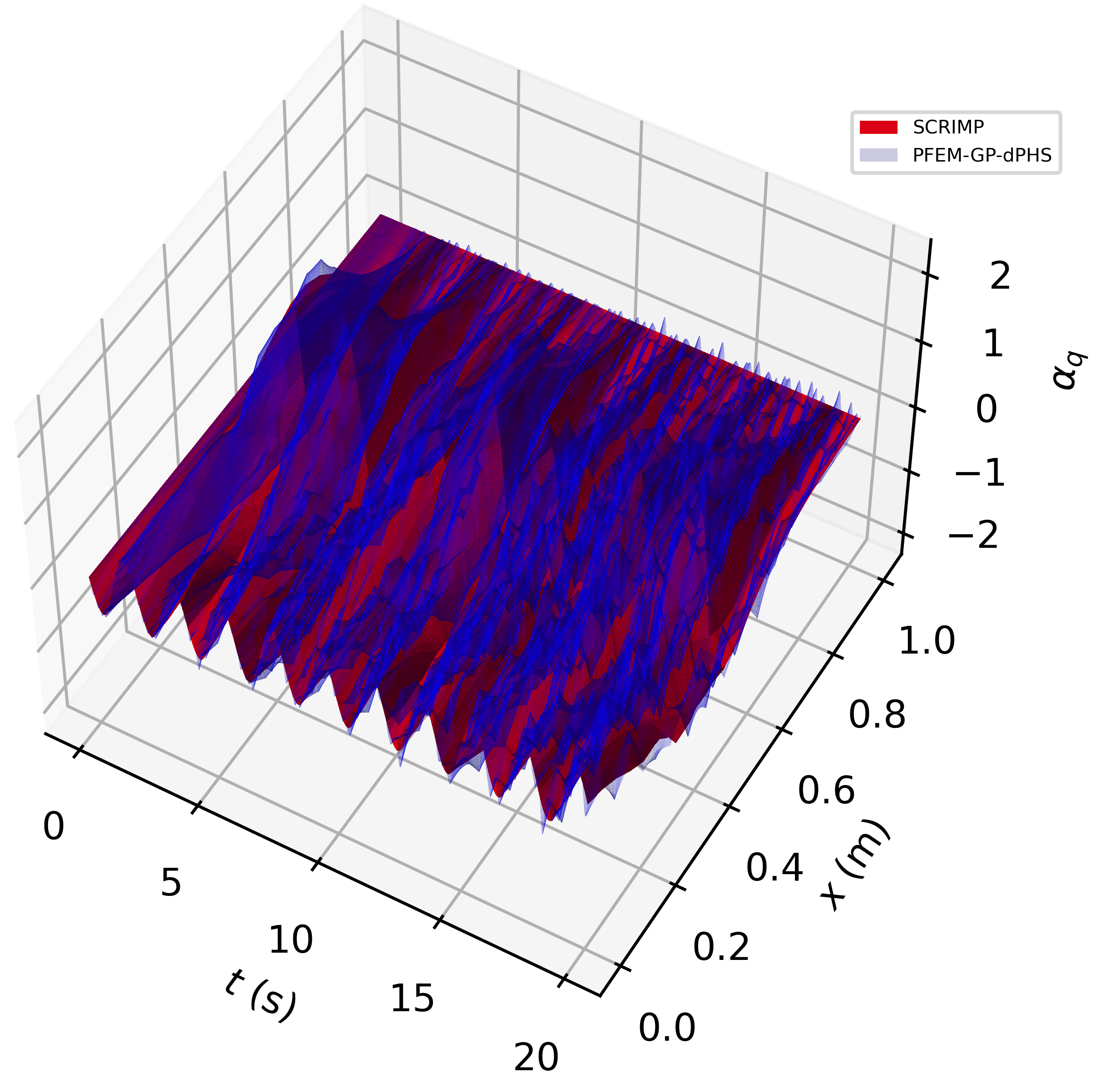}\centering
    \includegraphics[width=0.45\textwidth]{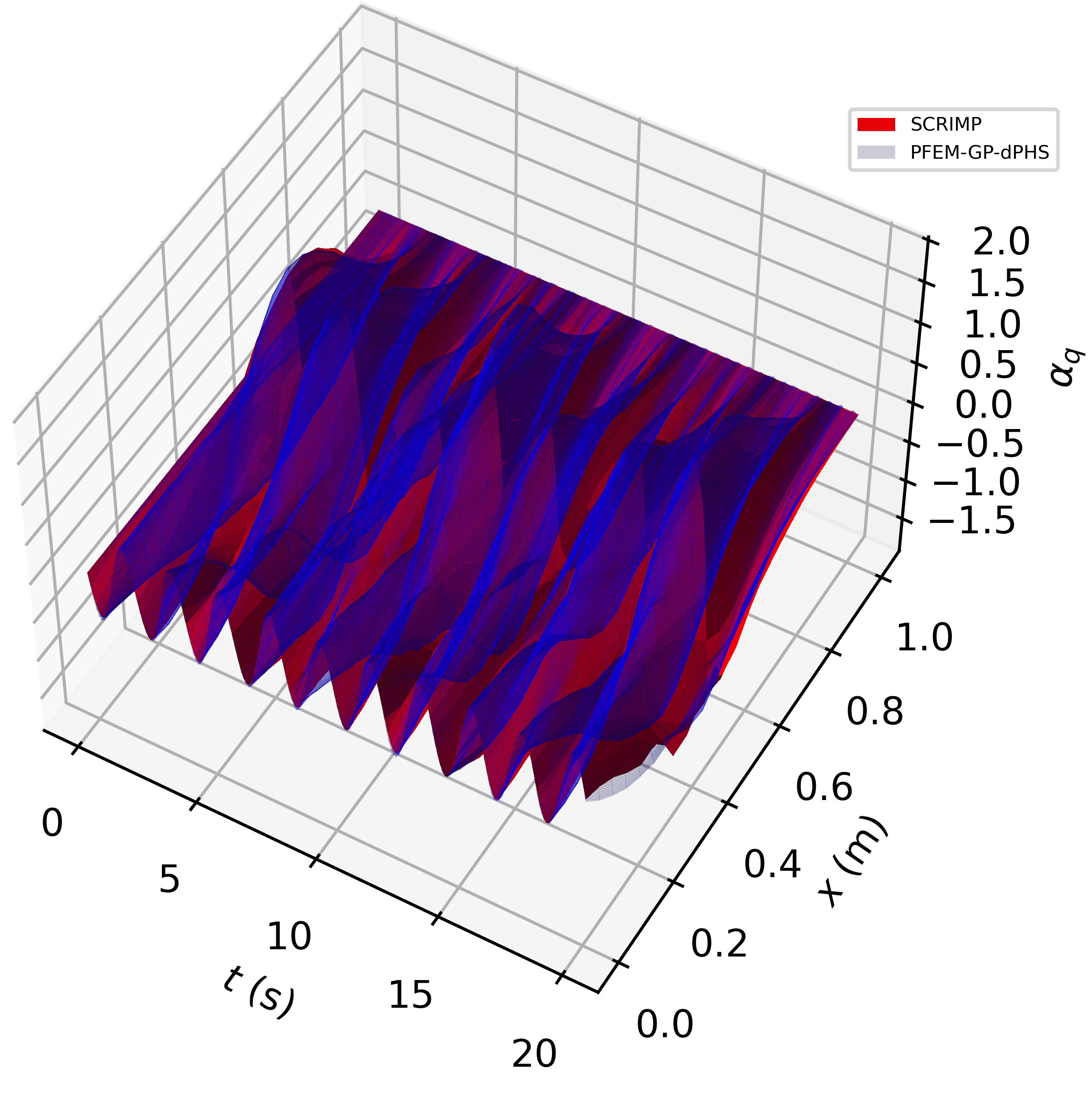}
    \caption{True $\alpha_q$ trajectory obtained with \texttt{SCRIMP} plotted with the estimated $\alpha_q$ trajectory obtained using the posterior mean of the PFEM-GP-dPHS and discretizing hyperparameters. Left: a result using the same parameters and hyperparameters as those of Figure \ref{fig:hyperparameters} ($N_q = N_p = 21$ and $\Delta x = 0.2$). Right:  a result using finer parameters, $N_q = N_p = 41$ and $\Delta x = 0.1$ resulting in 46 hyperparameters ($N_s = 35$ for both). This shows that using a larger number of finite elements (right) for $\alpha_q$ and $\alpha_p$ allows to obtain smoother simulated trajectories.}
    \label{fig:alpha_q_mix}
\end{figure}
\end{center}

\subsection{Computational costs}\label{subsec:computational_costs} 
The computational cost of the method is dominated by the inversion of the square kernel Gram matrix of size $(N_q+N_p)N_s$ using its Cholesky decomposition, for the calculation of the NLML resulting in cost of $O(((N_q+N_p)N_s)^3)$. Moreover, when calculating the gradient of the NLML, one needs one matrix multiplication of matrix of size $(N_q+N_p)N_s$ per hyperparameters of type $\ell_q$, $\ell_p$, $\sigma_f$ and $\sigma_{\text{noise}}$. This gives a gradient cost dominated by $O((N_{\theta q}+N_{\theta p}+2)((N_q+N_p)N_s)^3)$ where $N_{\theta q}$ and $N_{\theta p}$ are the number of hyperparameters of type $\ell_q$ and $\ell_p$. This shows the advantage of using late lumping : using early lumping would imply $N_{\theta q} = N_q$ and $N_{\theta p} = N_p$ whereas the late lumping approach allows us to use lower dimensional discretizations of the functional hyperparameters. This reduces the number of hyperparameters to optimize, reduces the computational complexity of NLML optimization, and accelerates the whole computational procedure. Using an 8 year old personal laptop with a processor \verb|Intel Core i7-6700HQx8|, 20 minutes were needed to obtain the results of Figure \ref{fig:hyperparameters}, 1 hour for Figure \ref{fig:alpha_p_mix} (right), 10 minutes for the polynomial approximation of Figure \ref{fig:gp_polyn} and 2 hours and 30 minutes for the 2D case of Figure \ref{fig:2D_error}.

\section{Exploiting the posterior variance provided by PFEM-GP-dPHS}\label{app:posterior_variance}
Figure \ref{fig:exploit_variance} shows an example of the \emph{absolute} $L^2$ errors $e^{\text{abs}}_{i,j} = \|\alpha_p^i(t_j) - \alpha_p^{\text{SCP}}(t_j)\|_{L^2}$ for some example $i$ as a function of time, as opposed to Figure \ref{fig:rel_error_alpha_p} (left) which shows \emph{relative} errors (see Section \ref{sec:num} for further details on the definition of the quantities above). Figure \ref{fig:exploit_variance} also shows the quantities $\text{var}_q(t)$ and $\text{var}_p(t)$, which are given in equations \eqref{eq:var_q} and \eqref{eq:var_p} below :
\begin{align}
\text{var}_q(t) &\coloneqq \text{Tr}\Big(\Cov\big[\widetilde{\underline{e}_q}\big(\underline{\alpha}(t)\big)\big]\Big) = \mathbb{E}\big[\|\widetilde{\underline{e}_q}\big(\underline{\alpha}(t))-\mu_q(t)\|_2^2\big], \ \ \ \mu_q(t) = \mathbb{E}[\widetilde{\underline{e}_q}\big(\underline{\alpha}(t))]\,,\label{eq:var_q} \\
\text{var}_p(t) &\coloneqq \text{Tr}\Big(\Cov\big[\widetilde{\underline{e}_p}\big(\underline{\alpha}(t)\big)\big]\Big) = \mathbb{E}\big[\|\widetilde{\underline{e}_p}\big(\underline{\alpha}(t))-\mu_p(t)\|_2^2\big], \ \ \ \mu_p(t) = \mathbb{E}[\widetilde{\underline{e}_p}\big(\underline{\alpha}(t))]\,.\label{eq:var_p} 
\end{align}
 Above, $(\widetilde{\underline{e}_q}(\underline{\alpha}))_{\underline{\alpha}\in\R^{N_q+N_p}}$ is the GP
obtained by conditioning the GP $\underline{e}_q(\underline{\alpha})$ on the training data. Recall that $(\underline{e}(\underline{\alpha}))_{\underline{\alpha}\in\R^{N_q+N_p}}$ is the vector-valued GP with mean and kernel functions given in equations \eqref{eq:mean_e} and \eqref{eq:kernel_e} respectively, and $\underline{e}(\underline{\alpha}) = \big(\underline{e}_q(\underline{\alpha}), \underline{e}_p(\underline{\alpha})\big)$.
The equality $\text{Tr}(\Cov[\widetilde{\underline{e}_q}(\underline{\alpha}(t))]) = \mathbb{E}[\|\widetilde{\underline{e}_q}(\underline{\alpha}(t))-\mu_q(t)\|_2^2]$ (likewise for $\widetilde{\underline{e}_p}$) are justified at the end of this section.
\begin{figure}[t!]
    \centering
    \includegraphics[width=0.45\linewidth]{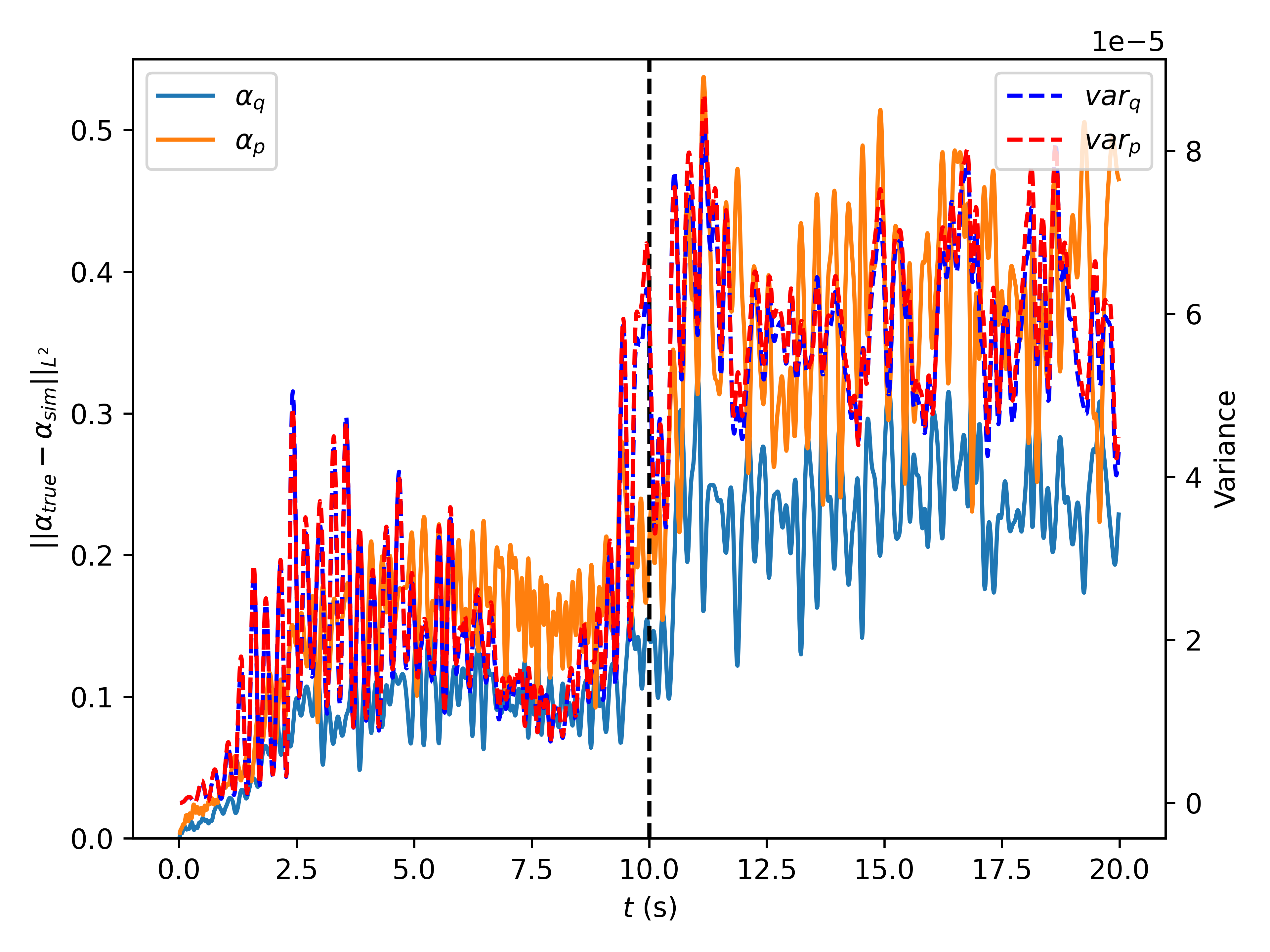}
    \caption{Evolution of $\text{var}_q(t),\ \text{var}_p(t)$ and the $L^2$ absolute errors as functions of time.}
    \label{fig:exploit_variance}
\end{figure}
Figure \ref{fig:exploit_variance} shows that the absolute error plots behave approximately like the quantities $\text{var}_q(t)$ and $\text{var}_p(t)$. In particular, all four curves increase after the 10s landmark, which corresponds to the time at which we stopped collecting measurement data in the \texttt{SCRIMP} simulation to train our PFEM-GP-dPHS model. 
The knowledge of this posterior variance can be used to perform adaptive learning for our PFEM-GP-dPHS model, as dexcribed in \cite[Section 5.4]{Courteville2000500} in a (finite-dimensional) GP-PHS context.

Note that in Figure \ref{fig:exploit_variance} shows error plots about the energy variables $(\underline{\alpha}_q,\underline{\alpha}_p)$, while the variance plots concern the \emph{coenergy} variables $(\underline{e}_q, \underline{e}_p)$ and not $(\underline{\alpha}_q,\underline{\alpha}_p)$. This is because we have no GP model over $\underline{\alpha} = (\underline{\alpha}_q,\underline{\alpha}_p)$; they are the inputs of the GP $(\underline{e}(\underline{\alpha}))_{\underline{\alpha}\in\R^{N_q+N_p}}$! However, after training the GP $\underline{e}(\underline{\alpha})$ and during a PFEM-GP-dPHS simulation, the variance of the conditioned GP $\widetilde{\underline{e}}$ \emph{can} be computed online at points $\underline{\alpha}(t)$, for any $t\geq 0$.

We finish this appendix with a proof of the equalities given in equations \eqref{eq:var_q} and \eqref{eq:var_p}. For this, consider a random vector $X\in\R^N$ with mean $m$ and covariance matrix $\Sigma$, then using that $\text{Tr}(AB) = \text{Tr}(BA)$ if $A$ and $B$ have suitable dimensions, 
\begin{align*}
    \text{Tr}\ \Sigma &= \text{Tr} \ \mathbb{E}[(X-m)(X-m)^\top] = \mathbb{E}\Big[ \text{Tr}\Big( (X-m)(X-m)^\top\Big)\Big] \\
    &=\mathbb{E}\Big[ \text{Tr}\Big( (X-m)^\top(X-m)\Big)\Big] = \mathbb{E}\big[ \text{Tr}\big( \|X-m\|_2^2\big)\big] =\mathbb{E}\big[\|X-m\|_2^2\big]\,.
 \end{align*}
To obtain the equalities in equations \eqref{eq:var_q} and \eqref{eq:var_p}, use the result above with $(N,X) = (N_q,\widetilde{\underline{e}_q}(\underline{\alpha}(t))$ and 
$(N,X) =  (N_p,\widetilde{\underline{e}_p}(\underline{\alpha}(t)))$.
\end{document}